\documentclass{birkau}

\usepackage{amsmath,amssymb,url}
\usepackage{enumitem} 
\usepackage{kantlipsum}
\usepackage{microtype}
\usepackage{bm}
\usepackage{tikz}
\usepackage{mathrsfs}
\usepackage{mathabx}
\usepackage{color}

\numberwithin{equation}{section}

\theoremstyle{plain}
\newtheorem{theorem}{Theorem}[section]
\newtheorem{lemma}[theorem]{Lemma}
\newtheorem{proposition}[theorem]{Proposition}
\newtheorem{corollary}[theorem]{Corollary}
\newtheorem{example}[theorem]{Example}

\theoremstyle{definition}
\newtheorem{definition}[theorem]{Definition}

\newtheorem{remark}[theorem]{Remark}
 
\DeclareMathSizes{10}{10}{6}{6}

\newcommand{\e}{\mathsf e}

\begin{document}

\title[Universal clone algebra]{Universal clone algebra} 
\author[A. Salibra]{Antonino Salibra}
\address{Institut de Recherche en Informatique Fondamentale\\
Universit\'e Paris Cit\'e\\ 8 Place Aur\'elie Nemours, 75205 Paris Cedex 13, France}
\email{salibra@unive.it}
\urladdr{http://www.dsi.unive.it/~salibra}

\subjclass{Primary: 08A40; Secondary: 08B05, 08B15, 08C05}

\keywords{Universal clone algebra, t-Algebras, Clone algebras, Functional clone algebras, Universal algebra, Clones, Infinitary algebras, Birkhoff's theorem, Topological Birkhoff's theorem}

\begin{abstract} 
 We develop a new general framework for algebras and clones, called universal clone algebra. Algebras and clones of finitary operations are to universal algebra what t-algebras and clone algebras are to universal clone algebra. Clone algebras have been recently introduced to found a one-sorted, purely algebraic theory of clones, while t-algebras are first introduced in this article.
We present a method to codify algebras and clones into t-algebras and clone algebras, respectively. We provide concrete examples showing that general results in universal clone algebra, when translated in terms of algebras and clones, give new versions of known theorems in universal algebra. We apply this methodology to  Birkhoff's HSP theorem and to the recent topological versions of Birkhoff's theorem.
\end{abstract}

\maketitle


\vspace{-1cm}

\section{Introduction}\label{sec:intro}

Clones are sets of finitary operations on a fixed carrier set that contain all projection operations and are closed under composition. They play an important role in universal algebra, since the set of all term operations of an algebra always forms a clone and in fact every clone is of this form. Moreover, important properties, like whether a given subset forms a subalgebra, or whether a given map is a homomorphism, depend on the clone of its term operations. Hence, comparing clones of algebras is much more appropriate than comparing their basic operations, in order to classify algebras according to different behaviours (see \cite{SZ86,T93}).

Clones play another important role in the study of first-order structures. 
Indeed, the polymorphism clone of a first-order structure, consisting of all finitary functions which preserve the structure, carry information about the structure, and is a powerful tool in its analysis.
Clones are also important in theoretical computer science. Many computational problems can be phrased as constraint satisfaction problems (CSPs): in such a problem, we fix a structure $A$ (also called the template or constraint language). The problem $\mathrm{CSP}(A)$ is the computational problem of deciding whether a given conjunction of atomic formulas over the signature of $A$ is satisfiable in $A$. The seminal discovery in the algebraic approach to CSP is the result of Jeavons \cite{J98} that, for a finite structure $A$, the complexity of $\mathrm{CSP}(A)$ is completely determined by the polymorphism clone of  $A$ (see e.g. \cite{BJK05,B15,BKW17}).

Some attempts  have been made to encode clones into algebras. A particularly important one led to  the concept of abstract clones  \cite{T93}, which are many-sorted algebras axiomatising composition of finitary functions and projections. Every abstract clone has a concrete representation as an isomorphic clone of finitary operations. 
Another important attempt to encode clones into algebras is due to Neumann \cite{neu70} and led to the concept of  abstract $\aleph_0$-clones, which are infinitary algebras axiomatising all projections and one infinitary operation of composition of arity $\omega$. Every abstract $\aleph_0$-clone has a concrete representation as an isomorphic $\aleph_0$-clone of infinitary term operations.



Bucciarelli and the author have recently introduced  in \cite{BS22} a one-sorted algebraic  theory of clones.  
Indeed, clone algebras ($\mathsf{CA}$) are defined by universally quantified equations and thus form a variety in the universal algebraic sense. 
 A crucial feature of this approach is connected to the role played by
variables in free algebras and  projections in clones. In clone algebras these are abstracted out, and take  the form of a countable infinite system of fundamental elements (nullary operations) $\e_1,\e_2,\dots,\e_n,\dots$ of the algebra. 
 One important consequence of the abstraction of variables and projections is the abstraction of
 term-for-variable substitution and  functional composition in $\mathsf{CA}$s, obtained  by introducing an $(n+1)$-ary operator $q_n$ for every $n\geq 0$. Roughly speaking, $q_n(a,b_1,\dots,b_n)$  represents
the substitution of $b_i$ for  $\e_i$ into $a$ for $1\leq i\leq n$ (or the composition of  $a$ with  $b_1,\dots,b_n$ in the first $n$ coordinates of $a$).

The most natural $\mathsf{CA}$s, the ones the axioms are intended to characterise, are algebras of functions, called \emph{functional clone algebras} 
($\mathsf{FCA}$s). 
 The elements of a $\mathsf{FCA}$ with value domain $A$ are infinitary operations, called here t-operations. They are functions $\varphi : \mathsf a\to A$,  whose domain $\mathsf a$, called a \emph{trace on $A$}, is a nonempty subset of $A^\mathbb N$ satisfying the following condition:
 \begin{equation}\label{eqtrace}\forall r,s\in A^\mathbb N.\ s\in \mathsf a\ \text{and}\ |\{i: s_i\neq r_i\}|<\omega  \Rightarrow r\in\mathsf a.\end{equation}
 A trace $\mathsf a$ on $A$ is \emph{complete} if $\mathsf a = A^\mathbb N$ and it is \emph{basic} if it is minimal with respect to $\subseteq$.
In this framework the nullary operators are the projections, defined by $\e_i^\mathsf a(s)=s_i$ for every $s\in \mathsf a$, and
 $q_n^\mathsf a(\varphi,\psi_1,\dots,\psi_n)$ represents the $n$-ary composition of $\varphi$ with $\psi_1,\dots,\psi_n$, acting on the first $n$ coordinates:
 $$q_n^\mathsf a(\varphi,\psi_1,\dots,\psi_n)(s)= \varphi(\psi_1(s),\dots,\psi_n(s),s_{n+1},s_{n+2},\dots),\ \text{for every $s\in \mathsf a$}.$$
Every clone algebra $\mathcal C=(C,q_n^\mathcal C,\e_i^\mathcal C)_{n\geq 0,i\geq 1}$  is isomorphic to a $\mathsf{FCA}$ with value domain $C$ and basic trace $\pmb \epsilon =\{s\in C^\mathbb N: |\{i: s_i\neq \e_i^\mathcal C\}|<\omega \}$.  

In \cite{BS22} the authors have shown that the finite-dimensional clone algebras ($\mathrm{Fd}\, \mathsf{CA}$) are the abstract counterpart of the clones of finitary operations, where the dimension of an element in a clone algebra is an abstraction of the notion of arity.
However, the most part of clone algebras are not finite-dimensional (e.g. the $\mathsf{FCA}$ of all t-operations). 
Then it is natural to investigate what are the algebraic structures that correspond to clone algebras in full generality (see the figure below).
\[
\left[
\begin{array}{ccc}
\mathrm{Fd}\,\mathsf{CA}  &  \text{Algebras} &\text{Clones}   \\
\mathsf{CA}  &  ? & ?  \\    
\end{array}
\right]
\]
In this paper  we introduce a new general framework for algebras and clones, called \emph{universal clone algebra}. To make a comparison, algebras and clones of finitary operations are to universal algebra what t-algebras and clone algebras are to universal clone algebra. Clone algebras have been already described in this introduction, while the notion of a t-algebra is first introduced  in this article.
A t-algebra is a tuple $\mathbf A=( A,\mathsf a,\sigma^\mathbf A)_{\sigma\in \tau}$, where $\tau$ is a set of operator symbols, $\mathsf a$ is a trace on the given set $A$ and $\sigma^\mathbf A:\mathsf a\to A$ is a t-operation for every $\sigma \in \tau$ (see the figure below). 
\[
\left[
\begin{array}{ccc}
\mathrm{Fd}\,\mathsf{CA}  &  \text{Algebras} &\text{Clones}   \\
\mathsf{CA}  & \text{t-Algebras} & \mathsf{FCA}s  \\    
\end{array}
\right]
\]
We have two algebraic levels. 
The lower degree of t-algebras and the higher degree of clone algebras. There are many ways to move between these levels, either individually or collectively.

If $\mathbf A$ is a t-algebra, then $\mathbf A^\uparrow$ is the term clone algebra over $\mathbf A$, the t-algebra analogue of the term clone of an algebra. 
If $\mathcal C$ is a clone $\tau$-algebra, then $\mathcal C^\downarrow$ is the t-algebra of type $\tau$ under $\mathcal C$, whose basic trace  is  generated by the sequence $(\e_1^\mathcal C,\dots,\e_n^\mathcal C,\dots)$ of the designated elements. We provide sufficient conditions in order that $\mathbf A^\updownarrows\cong \mathbf A$ and $\mathcal C^\downuparrows\cong \mathcal C$.

More interesting is the way we \emph{collectively} go up and down between the two levels (there is no relation with the  above described individual up-and-down).  If $K$ is a class of t-algebras of type $\tau$, then $K^\smalltriangleup$ is a class  of clone $\tau$-algebras. If $H$ is a class of clone $\tau$-algebras, we have two ways to go down: $H^\smalltriangledown$ and  $H^\blacktriangledown$ are two classes of t-algebras of type $\tau$ such that $H^\smalltriangledown\subseteq H^\blacktriangledown$.

We generalise the usual algebraic construction to t-algebras (namely, t-subalgebra, t-product, t-homomorphic image and t-variety) and define a class of t-algebras to be an Et-variety (resp. Ft-variety) if it is a t-variety closed under t-expansion (resp. under full t-expansion).
We prove that
\begin{itemize}
\item  If $K$ is a t-variety of t-algebras, then $K^\smalltriangleup$ is a  variety of clone algebras.
\item  If $H$ is a variety of clone algebras, then  $H^\smalltriangledown$ is an Ft-variety and  $H^\blacktriangledown$ is an Et-variety.
\end{itemize}
The main results of this article depend on  
the above described up-and-down technique  and
on the way we embed algebras into t-algebras. We provide concrete examples that general results in universal clone algebra, when translated in terms of algebras and clones, give new versions of known theorems in universal algebra. We apply this methodology to Birkhoff's HSP theorem and to the recent topological versions of Birkhoff's theorem.
We now describe the main results of this paper in more detail.

As well-explained in \cite{GP18}, there are two formulations of Birkhoff's theorem. The global formulation states that a class of algebras is a variety (i.e., closed under arbitrary products, subalgebras, and homomorphic images) if and only if it is an equational class (see \cite[Theorem 11.9]{BS81}).
The  local formulation of Birkhoff's theorem concerns the variety generated by a single algebra $\mathbf S$, i.e., the class of all algebras that can be obtained from $\mathbf S$ using arbitrary products, subalgebras, and homomorphic images. It states that an algebra $\mathbf T$ of the same type is a member of that variety if and only if it satisfies all equations that hold in $\mathbf S$, or in other words, if the natural clone homomorphism from the term clone $\mathrm{Clo}\mathbf S$ of $\mathbf S$ onto the term clone $\mathrm{Clo}\mathbf T$ of $\mathbf T$ exists (see \cite{Bir46}).

In this paper we first generalise the global version of Birkhoff's Theorem to Et-varieties and Ft-varieties. It takes the following form for a class $K$ of t-algebras of the same type: $K$ is an Et-variety (resp. Ft-variety) if and only if $K^\smalltriangleup$ is a variety of clone algebras and $K=K^{\smalltriangleup\blacktriangledown}$ (resp. $K=K^{\smalltriangleup\smalltriangledown}$).
As a corollary, we get the following  new version of Birkhoff HSP theorem for varieties of algebras.

\begin{theorem} {\rm (Birkhoff's Theorem for Algebras)}\label{thm:birkhoffforalgebras} Let $H$ be a class of algebras of the same type and $H^\star$ be the class of all t-algebras obtained by gluing together algebras in $H$ (formally defined in Section \ref{sec:bfafbfta}). Then the following conditions are equivalent:
\begin{enumerate}
\item $H$ is a variety of algebras; 
\item $H$ is an equational class of algebras;
\item  $H^\star$ is an Et-variety; 
\item $(H^\star)^\smalltriangleup$ is a variety of clone algebras and $H^\star=(H^\star)^{\smalltriangleup\blacktriangledown}$;
\item  $H^\star$ is an Ft-variety;
\item $(H^\star)^\smalltriangleup$ is a variety of clone algebras and $H^\star=(H^\star)^{\smalltriangleup\smalltriangledown}$;
\item  $H^\star$ is a t-variety.
\end{enumerate}
\end{theorem}


The last sections of the paper are devoted to the study of topological variants of Birkhoff's theorem, as initiated  by Bodirsky and Pinsker \cite{BP15} for locally oligomorphic algebras, and generalised recently by Schneider \cite{Sc17} and Gehrke-Pinsker \cite{GP18}.
These authors provide a Birkhoff-type characterisation of all those members $\mathbf T$ of the variety $\mathbb{HSP}(\mathbf S)$ generated by a given algebra $\mathbf S$, for which the natural homomorphism from $\mathrm{Clo}\mathbf S$  onto $\mathrm{Clo}\mathbf T$  is uniformly continuous with respect to the uniformity of pointwise convergence. 
Schneider \cite{Sc17} and Gehrke-Pinsker \cite{GP18} have independently   shown  that for any algebra $\mathbf T$ in the variety generated by an algebra $\mathbf S$, the induced natural clone homomorphism is uniformly continuous if and only if every finitely generated subalgebra of $\mathbf T$ is a homomorphic image of a subalgebra of a finite power of $\mathbf S$. 
The above theorem, termed Uniform Birkhoff in \cite{GP18}, has been successfully applied in the context of constraint satisfaction problems by Barto and Pinsker \cite{BP16}, and by Barto, Opr\v{s}al and Pinsker \cite{BOP18}.
We generalise this theorem to t-algebras as follows.
 Let $\mathbf A, \mathbf  B$ be t-algebras of the same type and let $\mathsf b$ be the trace of $\mathbf B$. Then the following
are equivalent.
\begin{itemize}
\item $\mathbf  B$ is an element of the  Et-variety  generated by  $\mathbf A$, and the natural homomorphism from the term clone algebra $\mathbf A^\uparrow$ onto the term clone algebra $\mathbf B^\uparrow$  is uniformly continuous.
\item Every t-subalgebra $\mathbf B_{\bar s}$ of $\mathbf  B$ generated by $s\in\mathsf b$  is a t-homomorphic image of a t-subalgebra of a finite t-power of $\mathbf A$.
\end{itemize}

%
%

We remark that the t-subalgebras involved in (2) depend on the trace $\mathsf b$. For example, if $s\in\mathsf b$ and  $|\{s_i : i\in\mathbb N\}|=\omega$, then $\mathbf B_{\bar s}$ is not in general finitely generated.
As a corollary, besides the version of topological Birkhoff by Schneider \cite{Sc17} and Gehrke-Pinsker \cite{GP18}, we get new versions of the topological Birkhoff's theorem for algebras depending on the choice of the trace $\mathsf b$.

We now describe the plan of this work.
In Sections \ref{sec:prelim}  we present some preliminary notions on algebras, clones, uniform spaces and Neumann's abstract $\aleph_0$-clones.
 In Section \ref{sec:tt}  we introduce the notions of thread, trace, t-operation and top extension of a finitary operation. In Section \ref{sec:clonealg} we recall  the basic concepts related to clone algebras.
Section \ref{sec:ta}  is devoted to the notion of a t-algebra. 
 The fundamental algebraic constructions of algebras are here generalised to t-algebras.
In Section \ref{sec:fca} functional clone algebras are defined in full generality.  
We introduce the t-operations determined by the elements of a clone algebra, via the operators $q_n$.
We use these t-operations to represent every clone algebra as a functional clone algebra of basic trace.  The t-algebra under a given clone algebra is also described here.
In Section \ref{sec:cats} the clone algebras of terms and hyperterms are introduced. The t-algebra analogue of the term clone of an algebra is also defined.
In Section \ref{sec:tvarieties} t-varieties, Ft-varieties and Et-varieties are established. We also describe the fundamental notion of a t-subalgebra generated by a thread. The up and down technique is described in Section \ref{sec:updown}. Birkhoff's theorems for t-algebras are proved in Sections \ref{sec:gb} and \ref{sec:gb2}, while in Section \ref{sec:newbir} a technique for building up t-algebras from algebras is described. A new Birkhoff's theorem for variety of algebras is a consequence of this technique.
The last two sections of the paper are devoted to free t-algebras and to a topological version of Birkhoff's theorem for t-algebras.


\section{Preliminaries}\label{sec:prelim}

The notation and terminology in this paper are pretty standard. For
concepts, notations and results not covered hereafter, the reader is
referred to \cite{BS81,mac87} for universal algebra, to \cite{L06,SZ86,T93} for the theory of clones and to \cite{J87} for uniform spaces.

In this paper $\mathbb N=\{1,2,\dots\}$ denotes the set of positive natural numbers.

By a \emph{finitary operation} on a set $A$ we  mean  a function $f:A^n\to A$ for some natural number $n\geq 0$. 
We denote by $O_A$  the set of all finitary operations on a set $A$. 
If $F\subseteq O_A $, then $F^{(n)}= \{f:A^n\to A\ |\ f\in  F\}$. 

 By an \emph{infinitary operation} on $A$ we always mean a function from $A^\mathbb N$ into $A$. The set of all infinitary operations on $A$ is denoted by $A^{A^\mathbb N}$.
 
 If $f:A\to B$ is a function, then we denote by $f^\mathbb N$ the function from $A^\mathbb N$ into $B^\mathbb N$ such that $f^\mathbb N(s)= (f(s_i) : i\in\mathbb N)$ for every $s\in A^\mathbb N$.
 
 If $\theta$ is an equivalence relation on a set $A$, then $[a]_\theta$ denotes the equivalence class of $a\in A$.

\subsection{Algebras}\label{sec:alg}
A finitary type is an $\omega$-sorted set $\rho=(\rho_n:n\in \omega)$ of disjoint sets of operator symbols. 
A  $\rho$-algebra is a tuple $\mathbf S=(S,\sigma^\mathbf S)_{\sigma\in\rho}$ such that $\sigma^\mathbf S:S^n\to S$ is a finitary operation of arity $n$ for every $\sigma\in\rho_n$. As a matter of notation, algebras will be denoted by the letters $\mathbf R, \mathbf S,\dots$.

In the following we fix a countable infinite set $I=\{v_1, v_2,\dots,  v_n,\dots\}$ of variables.
We denote by $F_\rho(I)$ the set of  $\rho$-terms over the countable infinite set $I$ of variables.
If $t$ is a $\rho$-term, then we write $t=t(v_1,\dots,v_n)$ if $t$ can be built up starting from variables $v_1,\dots,v_n$.

If $K$ is a variety of $\rho$-algebras, then we denote by $\mathbf F_K$ its free algebra over the countable infinite set
$I$ of generators. 

If $K$ is a class of $\rho$-algebras, then $\mathrm{Th}_K$ denotes the set of identities between $\rho$-terms true in every member of $K$.
If $K=\{\mathbf S\}$ is a singleton, then we write $\mathrm{Th}_\mathbf S$ for $\mathrm{Th}_{\{\mathbf S\}}$.

Closure of a class of similar algebras under  homomorphic images, direct products, subalgebras and  isomorphic images is denoted by $\mathbb{H}$, $\mathbb{P}$, $\mathbb{S}$ and $\mathbb{I}$ respectively. Closure under finite products is denoted by $\mathbb{P}^{\mathrm{fin}}$.
$\mathrm{Var}(K)$ denotes the variety generated by a class $K$ of similar algebras. 


{\bf Warning}: In Section \ref{sec:ta} of this article we introduce new algebraic structures, called t-algebras. In the following the words algebra and $\rho$-algebra  always mean algebra of finitary operations and algebra of a finitary type $\rho=(\rho_n: n\geq 0)$. 
If there is danger of confusion, we write finitary algebra for algebra. The reader should be aware that  ``finite algebra'' has another meaning: it is an algebra of finite cardinality.

\subsection{Clones of operations}\label{sec:clo} 

\bigskip

We recall some notations and terminology on clones we will use in the following. 
 The \emph{composition}  of $f\in  O_A^{(n)}$ with  $g_1,\dots,g_n\in   O_A^{(k)}$ is the operation $f(g_1,\dots,g_n)\in  O_A^{(k)}$ defined as follows, for all $\mathbf a\in A^k$: 
$$f(g_1,\dots,g_n)(\mathbf a)= f(g_1(\mathbf a),\dots, g_n(\mathbf a)).$$

  A \emph{clone on a set $A$} is a subset $F$ of $O_A$ containing all projections
$p^{(n)}_i:A^n\to A$ ($n\geq i$) and closed under composition. 
A \emph{clone on a $\rho$-algebra $\mathbf S$} is a clone on $S$  containing  the operations $\sigma^\mathbf S$ ($\sigma\in\rho$) of $\mathbf S$.
 
We denote by $\mathrm{Clo}\mathbf S$ the term clone of an algebra $\mathbf S$.

 The classical approach to clones, as evidenced by the standard monograph \cite{SZ86}, considers clones only containing operations that are at least unary. 
However, the full generality of some results in this paper requires clones allowing nullary operators.

Let $C$ and $D$ be clones. 
A \emph{clone homomorphism} from $C$ into $D$ is a mapping $F : C\to D$ which
\begin{itemize}
\item sends functions of $C$ to functions of the same arity in $D$,
\item sends projections in $C$ to the corresponding projections in $D$, and
\item preserves composition, i.e., for all $f, g_1,\dots, g_n \in C$, $F(f(g_1,\dots,g_n)) = F(f)(F(g_1),\dots,F(g_n))$.
\end{itemize}

\subsubsection{Neumann's abstract $\aleph_0$-clones \cite{neu70,T93}}\label{sec:neu} 
The idea here is to regard an $n$-ary operation $f$ as an infinitary operation that only depends on the first $n$ arguments (see Definition \ref{def:topext}). The corresponding abstract definition is as follows.
An \emph{abstract $\aleph_0$-clone} is an infinitary algebra $(A,\e_i,q_\infty)$, where the $\e_i$ ($1\leq i<\omega$) are nullary operators and $q_\infty$ is an infinitary operation of arity $\omega$ satisfying the following axioms: 
\begin{enumerate}[label=\textup{N\arabic*.},leftmargin=1.25\parindent]
\item 
$q_\infty(\e_i,x_1,\dots,x_n,\dots)=x_i$;
\item  
$q_\infty(x,\e_1,\dots,\e_n,\dots)=x$;
\item 
$q_\infty(q_\infty(x,\mathbf  y),\mathbf z)= q_\infty(x,q_\infty(y_1,\mathbf z),\dots,q_\infty(y_n,\mathbf z),\dots)$,
where  $\mathbf y$ and $\mathbf  z$ are countable infinite sequences of variables.  
\end{enumerate}

A \emph{functional $\aleph_0$-clone} with value domain $A$ is an algebra $(F,\e_i^{A^\mathbb N},q_\infty^{A^\mathbb N})$, where $F\subseteq  A^{A^\mathbb N}$ is a set of infinitary operations on $A$ and, for every $\varphi,\psi_i\in F$ and $s\in A^\mathbb N$, $\e_i^{A^\mathbb N}(s)=s_i$ and 
$q_\infty^{A^\mathbb N}(\varphi,\psi_1,\dots,\psi_n,\dots)(s)=\varphi(\psi_1(s),\dots,\psi_n(s),\dots)$. 

Neumann shows in \cite{neu70} that every abstract $\aleph_0$-clone is isomorphic to a functional $\aleph_0$-clone, which is the free algebra over the set $\{\e_1,\dots,\e_n,\dots\}$ of generators in a variety of infinitary algebras defined by operations of arity $\omega$.

There is a faithful functor from the category of clones to the category of abstract $\aleph_0$-clones, but this functor is not onto. 

The connection between Neumann's abstract $\aleph_0$-clones and clone algebras was explained in \cite[Section 4.3]{BS22}.

\subsection{Uniform spaces}\label{sec:us}
A \emph{uniformity} on a set $X$ is a filter $U$ of reflexive binary relations on $X$  satisfying the following conditions:
\begin{itemize}
\item If $\alpha$ is a member of $U$, then so is $\alpha^{-1} =\{(y,x) : (x,y)\in\alpha\}$.
\item For every $\alpha\in U$, there exists some $\beta\in U$ such that $\beta\circ\beta\subseteq \alpha$, where $\circ$ is the composition of binary relations.
\end{itemize}
The elements of $U$ are called \emph{entourages}. If $\alpha$ is an entourage and $x\in X$, then $\alpha[x] := \{y\in  X : (x, y) \in \alpha\}$.
For a subset $Y$ of $X$, $\alpha[Y]= \bigcup_{x\in Y} \alpha[x]$.


 For a uniformity $U$ on $X$, the \emph{uniform topology on $X$} is the topology defined as follows: $S$ is open if, for every $x\in S$, there exists an entourage $\alpha$ such that $\alpha[x] \subseteq S$. 

  The \emph{closure} of a set $T$ is the set $\overline T=\cap_{\alpha\in U}\alpha[T]$.
 Every set $A$ containing $\alpha[x]$ for some entourage $\alpha$ is a neighbourhood of $x$. In particular, the family $\alpha[x]$ ($\alpha\in U$) is a fundamental system of neighbourhoods of $x$.

 Let $X, Y$ be uniform spaces. A map $f : X \to Y$ is called \emph{uniformly continuous} if for every entourage $\alpha$ of $Y$ there exists some entourage $\beta$ of $X$ such that $(f \times f)(\beta)\subseteq \alpha$.

\section{Threads and Traces}\label{sec:tt}
In this section we introduce the notion of a \emph{trace}, which is necessary to define  functional clone algebras and t-algebras.

Let $A$ be a set. Every element $r\in A^\mathbb N$ is called a \emph{thread on $A$}.
If $r$ is a thread, then $\mathrm{set}(r)=\{r_i : i\in\mathbb N\}$.

As a matter of notation, if $r$ is a thread on $A$ and $a_1,\dots,a_n\in A$, then $r[a_1,\dots,a_n]$ is a thread defined as follows:
 $$r[a_1,\dots,a_n](i)=\begin{cases}a_i&\text{if $i\leq n$}\\ r_i&\text{if $i > n$}\end{cases}$$

\subsection{Traces} We define an equivalence relation $\equiv_\mathbb N$ on the set $A^\mathbb N$ of all threads as follows:
$$r \equiv_\mathbb N s\ \text{iff}\ |\{i: r_i\neq s_i\}| < \omega.$$
In other words, $r \equiv_\mathbb N s$ iff there exist $n\in\mathbb N$ and  $a_1,\dots,a_n\in A$ such that $s=r[a_1,\dots,a_n]$.
$[r]_\mathbb N^A$ denotes the equivalence class of the thread $r\in A^\mathbb N$. If there is no danger of confusion, we write $[r]_\mathbb N$ for $[r]_\mathbb N^A$.

\begin{definition}\label{def:3.1}
A subset $\mathsf a$ of $A^\mathbb N$ is called a \emph{trace on $A$} if it is union of equivalence classes of $\equiv_\mathbb N$:
$$s\in \mathsf a\ \text{and}\ r\equiv_\mathbb N s \Rightarrow r\in\mathsf a.$$
\end{definition}
\noindent A trace $\mathsf a$ on $A$ is called  \emph{basic}  if $\mathsf a=[s]_\mathbb N$ for some $s\in A^\mathbb N$;
 \emph{complete}  if $\mathsf a=A^\mathbb N$.

If $\mathsf a$ is a trace on $A$, then  $\mathsf d\subseteq_{\mathrm{bas}} \mathsf a$ means that $\mathsf d$ is a basic trace and $\mathsf d\subseteq \mathsf a$.

\begin{definition}\label{def:t-op} Let $A$ be a set and $\mathsf a$ be a trace on $A$.
 A \emph{t-operation}  is a map $\varphi: \mathsf a\to A$.
\end{definition}

 We denote by $A^{\mathsf a}$ the set of all t-operations from a trace $\mathsf a$ into $A$. 
 
 If $b\in A$, then $c_b^\mathsf a:\mathsf a\to A$ denotes the constant function such that $c_b^\mathsf a(s)=b$ for every $s\in\mathsf a$.
 
\begin{definition}\label{def:semi}
  A t-operation $\varphi: \mathsf a\to A$ is \emph{semiconstant} if, for every $r,s\in\mathsf a$,  $r\equiv_\mathbb N s\Rightarrow\varphi(s)=\varphi(r)$.
 \end{definition}

\subsection{Similarity and blocks}
\begin{definition}\label{def:topext}
 Let $\mathsf a$ be a trace on a given set $A$ and $f:A^n\to A$ be a finitary operation. The  \emph{top $\mathsf a$-extension $f_\mathsf a^\top:\mathsf a\to A$ of $f$} is defined as follows:
 $$f_\mathsf a^\top(s)=f(s_1,\dots,s_n),\ \text{for every $s\in \mathsf a$}.$$
 We say that $f,g\in  O_A$ are \emph{similar}, and we write $f\approx_{O_A} g$, if   $f_\mathsf a^\top=g_\mathsf a^\top$ for some (and then all) trace $\mathsf a$ on $A$.
\end{definition}

  The set of all top $\mathsf a$-extensions is denoted by $B^\top_\mathsf a=\{f_\mathsf a^\top : f\in O_A\}$.



 If $B$ is a block (i.e., equivalence class) of $\approx_{O_A}$, then $B \cap  O^{(n)}_A$ is either empty or a singleton.  
If $B \cap  O^{(n)}_A\neq\emptyset$, then $B \cap  O^{(k)}_A\neq\emptyset$ for every $k\geq n$. 


\begin{lemma}  \cite[Lemma 3.13]{BS22}
  Every clone on $A$ is a union of blocks. 
\end{lemma}


\section{Clone algebras}\label{sec:clonealg}

In this section we recall from \cite{BS22} the definition of a \emph{clone algebra} ($\mathsf{CA}$) as a more canonical algebraic account of clones using standard one-sorted algebras. 
In our approach we replace Neumann's infinitary operator of composition by a countable infinite set of finitary operators of composition.

The finitary type of  clone algebras contains a countable infinite family of nullary operators $\e_i$ ($i\geq 1$) and, for each $n\geq 0$, an operator $q_n$ of arity $n+1$.  If $\tau$ is a set of operator symbols, then the algebraic type of clone $\tau$-algebras is $\tau \cup\{ q_n: n\geq 0\} \cup\{ \e_i : i\geq 1\}$, where each $\sigma\in\tau$ has arity $0$. 

In the remaining part of this paper 
when we write $q_n(x,\mathbf y)$ it will be implicitly stated that  $\mathbf y=y_1,\dots,y_n$ is a sequence of length $n$.

 \begin{definition} \label{def:clonealg}
 A \emph{clone $\tau$-algebra} ($\mathsf{CA}_\tau$, for short)  is an algebra 
 $$\mathcal C = (C, q^\mathcal C_n,\e^\mathcal C_i,\sigma^\mathcal C)_{n\geq 0,i\geq 1,\sigma\in\tau},$$ where $C$ is a set, $\sigma^\mathcal C$ and $\e^\mathcal C_i$ are elements of $C$, and $q^\mathcal C_n$ is an $(n+1)$-ary finitary operation satisfying the following identities:
\begin{enumerate}
\item[(C1)] $q_n(\e_i,x_1,\dots,x_n)=x_i$ $(1\leq i\leq n)$;
\item[(C2)] $q_n(\e_j,x_1,\dots,x_n)=\e_j$ $(j>n)$;
\item[(C3)] $q_n(x,\e_1,\dots,\e_n)=x$ $(n\geq 0)$;
\item[(C4)] $q_n(x,  \mathbf{y})= q_k(x, \mathbf{y},\e_{n+1},\dots,\e_k)$ ($k> n$);
 \item[(C5)] $q_n(q_n(x, y_1,\dots,y_n), \mathbf{z})=q_n(x,q_n(y_1, \mathbf{z}),\dots,q_n(y_n,\mathbf{z}))$.
\end{enumerate}
If $\tau$ is empty, an algebra satisfying (C1)-(C5) is called a \emph{pure clone algebra}.
\end{definition}

In the following, when there is no danger of confusion,
we will write $\mathcal C = (C, q^\mathcal C_n,\e^\mathcal C_i,\sigma^\mathcal C)$  for $\mathcal C = (\mathcal C, q^\mathcal C_n,\e^\mathcal C_i,\sigma^\mathcal C)_{n\geq 0,i\geq 1,\sigma\in\tau}$.

The class of clone $\tau$-algebras  is denoted by $\mathsf{CA}_\tau$ and the class of all clone  algebras of any type by $\mathsf{CA}$. $\mathsf{CA}_0$ denotes the class of all pure clone algebras.
We also use $\mathsf{CA}_\tau$ as shorthand for the phrase ``clone $\tau$-algebra'', and similarly for $\mathsf{CA}$. 
$\mathsf{CA}_\tau$ is a variety of algebras.


\begin{example}\label{exa:proj} The algebra $\mathcal P=(\mathbb N, q_n^\mathcal P,\e_i^\mathcal P)$, where  $\e_i^\mathcal P=i$ and $$q_n^\mathcal P(i,k_1,\dots,k_n)=\begin{cases}k_i&\text{if $i\leq n$}\\ i&\text{if $i>n$}\end{cases}$$ 
is the minimal pure  clone algebra.
\end{example}

\begin{example}\label{exa:free} Let $\rho$ be a finitary type, $K$ be a variety of $\rho$-algebras and  $\mathbf F_K=(F_K,\sigma^\mathbf F)_{\sigma\in\rho}$ be the free $K$-algebra over the countable set $I=\{v_1,v_2,\dots\}$ of generators. 
The set of all $n$-finite endomorphisms $f$ of $\mathbf{F}_K$ (i.e., $f(v_i)=v_i$ for every $i>n$) can be collectively expressed by an $(n+1)$-ary operation $q_n^\mathcal F$ on $F_K$ (see   \cite[Definition 3.2]{mac83} and \cite[Definition 5.2]{BS22}): $q_n^\mathcal F(a,b_1,\dots,b_n)=s(a)$ for every $a,b_1,\dots,b_n\in F_K$,
where $s$ is the unique $n$-finite endomorphism
of $\mathbf{F}_K$ which sends the generator $ v_i$ to $b_i$ ($1\leq i\leq n$). 
 Then the algebra  $\mathcal F_K=(F_K, q_n^\mathcal{F}, \e_i^\mathcal{F}, \sigma^\mathcal F)_{\sigma\in\rho}$ is a clone algebra,
where $\e_i^\mathcal{F}= v_i\in I$,  $q_n^\mathcal{F}$ is above defined and $\sigma^\mathcal F$ is the equivalence class in $\mathbf F_K$ of the  $\rho$-term $\sigma(v_1,\dots,v_n)$ ($\sigma\in\rho_n$).
\end{example}

\subsection{Independence and dimension}
We define the notions of independence and dimension in clone algebras, abstracting the notion of arity of the finitary operations. We follows \cite[Section 3]{BS22}.

\begin{definition} \cite[Definition 3.4]{BS22} An element $a$ of a clone algebra $\mathcal  C$ \emph{is independent of} $\e_{n}$ if 
$q_{n}^\mathcal C(a,\e_1^\mathcal C,\dots,\e_{n-1}^\mathcal C,\e_{n+1}^\mathcal C)=a$. If $a$ is not independent of $\e_{n}$, then we say that 
$a$ \emph{is dependent on} $\e_{n}$.
\end{definition}

By \cite[Lemma 3.5]{BS22} $a$ is independent of $\e_n$ iff $q_{n}^\mathcal C(a,\e_1^\mathcal C,\dots,\e_{n-1}^\mathcal C,b)=a$ for every $b\in C$. We define the \emph{dimension of $a$} as follows:
$$\mathrm{dim}(a)=\begin{cases}n&\text{if $n=\mathrm{max}\{k: \text{$a$ depends on $\e_k$}\}$} \\\omega&\text{if $\{k: \text{$a$ depends on $\e_k$}\}$ is infinite} \end{cases}
$$

An element $a\in C$ is said to be: (i) \emph{$k$-dimensional} if $\mathrm{dim}(a)=k$; (ii) \emph{finite-dimensional} if it is $k$-dimensional for some $k<\omega$; (iii) \emph{zero-dimensional}  if $\mathrm{dim}(a)=0$. If $a$ is zero-dimensional, then
$q_n^\mathcal C(a,b_1,\dots,b_n)=a$ for all $n$ and $b_1,\dots,b_n\in C$. We denote by $\mathrm{Zd}\mathcal C$ the set of all zero-dimensional elements of a clone algebra $\mathcal C$.

We denote by $\mathrm{Fd}\,\mathcal C$ the subalgebra of a clone algebra $\mathcal C$ constituted by  all its finite-dimensional elements. We say that $\mathcal C$ is \emph{finite-dimensional} if $\mathcal C=\mathrm{Fd}\,\mathcal C$.

The clone algebras introduced in Examples \ref{exa:proj} and \ref{exa:free} are all finite-dimensional.

\section{t-Algebras}\label{sec:ta}
In  \cite[Theorem 3.20]{BS22} the authors have shown that the finite-dimensional clone algebras are the abstract counterpart of the clones of finitary operations. This correspondence has been shown fruitful to characterise the lattices of equational theories of algebras and to study the category of all varieties of algebras.  Since the most part of clone algebras are not finite-dimensional, then it is natural to determine the algebraic structures that correspond to clone algebras in full generality. This is done in this section by introducing t-algebras.

\subsection{The definition of a t-algebra} \label{sec:extfinalg}
The t-algebras are infinitary algebras. 

\begin{definition}\label{def:ddd} Let $\tau$ be a set of operator symbols. A tuple $\mathbf A= ( A,\mathsf a,\sigma^\mathbf A)_{\sigma\in \tau}$, where $A$ is a set, $\mathsf a$ is a trace on $A$ and $\sigma^\mathbf A: \mathsf a\to A$ is a t-operation for every $\sigma\in\tau$, is called a \emph{t-algebra of type $\tau$ and trace $\mathsf a$}. 
\end{definition}
As a matter of notation, t-algebras will be denoted by the letters $\mathbf A, \dots,  \mathbf G$.

A \emph{basic  (resp. complete) t-algebra} is a t-algebra whose trace is basic (resp. complete). 

\begin{example}\label{exa:diminf} Let $A=\{0,1\}$  and  $\bar 1$ be the sequence defined by $\bar 1(i)=1$ for every $i\in\mathbb N$. Then the tuple $\mathbf A=(A,[\bar 1]_\mathbb N,\sigma^\mathbf A)$, where for every $s\in  [\bar 1]_\mathbb N$ 
$$\sigma^\mathbf A(s)= \begin{cases}1 &\text{if $|\{s_i: s_i=0\}|$ is even}\\
0&\text{otherwise}\end{cases}$$
 is a basic t-algebra. Note that $\sigma^\mathbf A$ is not the top  extension of any finitary operation.
\end{example}

\begin{example}\label{exa:extfinalg} 
 Let $\rho$ be a finitary type,  $\mathbf S=(S,\sigma^\mathbf S)_{\sigma\in\rho}$ be a $\rho$-algebra and $\mathsf a$ be a trace on $S$. Then the tuple $\mathbf S_\mathsf a^\top=(S,\mathsf a,(\sigma^\mathbf S)_\mathsf a^\top)_{\sigma\in\rho}$ is called \emph{the t-algebra over $\mathbf S$ of trace $\mathsf a$ and type $\rho^\star=\bigcup_{n\geq 0} \rho_n$}.
\end{example}


In the remaining part of this section we extend to t-algebras the classical notions of subalgebra, homomorphism, congruence and product.

\subsection{t-Subalgebras} \label{sec:restr}
\begin{definition}
  Let $\mathbf A$ be a t-algebra of type $\tau$ and trace $\mathsf a$. A pair $(B,\mathsf b)$ is a \emph{t-subalgebra} of $\mathbf A$ if $B\subseteq A$, $\mathsf b\subseteq  \mathsf a$ is a trace on $B$ and $\sigma^\mathbf A(s)\in B$, for every $s\in \mathsf b$ and $\sigma\in\tau$.
\end{definition}

We denote by $\mathbf B=(B,\mathsf b, \sigma^\mathbf B)_{\sigma\in\tau}$, where  $\sigma^\mathbf B= (\sigma^\mathbf A)_{|\mathsf b}$, the t-subalgebra determined by the pair $(B,\mathsf b)$.  We write $\mathbf B\leq_t \mathbf A$ for ``$\mathbf B$ is a t-subalgebra of $\mathbf A$''. 

\begin{definition}\label{def:fullsub}
 A t-subalgebra $\mathbf B$ of $\mathbf A$ is called \emph{full} if $B=A$. If $\mathsf b$ is the trace of $\mathbf B$, then $\mathbf B$ will be denoted by 
$\mathbf A_{\upharpoonright\mathsf b}$. 
\end{definition}

If $K$ is a class of t-algebras of the same type, then we denote by $\mathrm S_t K$ the closure of $K$ by t-subalgebra.

\subsection{t-Homomorphisms} 
\begin{definition}
  Let $\mathbf A= ( A,\mathsf a,\sigma^\mathbf A)_{\sigma\in \tau}$  and $\mathbf B=( B,\mathsf b,\sigma^\mathbf B)_{\sigma\in \tau}$ be t-algebras of the same type.
  A map $f:A\to B$ is a \emph{t-homomorphism} from $\mathbf A$ into $\mathbf B$ if $f^\mathbb N: \mathsf a  \to \mathsf b$ and
$$f\circ \sigma^\mathbf A = \sigma^\mathbf B\circ f^\mathbb N,\ \text{for every $\sigma\in\tau$}.$$ 
\end{definition}
The \emph{image of $f$} is a t-subalgebra $\mathbf C=(C,\mathsf c, \sigma^\mathbf C)_{\sigma\in\tau}$ of $\mathbf B$ defined as follows: $C=\{ f(a): a\in A\}$, $\mathsf c= \{f^\mathbb N(s): s\in\mathsf a\}$ and $\sigma^\mathbf C(f^\mathbb N(s))= f(\sigma^\mathbf A(s))$.

Note that $f^\mathbb N$ respects the equivalence $\equiv_\mathbb N$ introduced in Definition \ref{def:3.1}: $r\equiv_\mathbb N s\Rightarrow f^\mathbb N(r)\equiv_\mathbb N f^\mathbb N(s)$. 

We say that a t-homomorphism $f$ from $\mathbf A$ into $\mathbf B$ is \emph{onto}
(and that $\mathbf B$ is a \emph{t-homomorphic image}  of $\mathbf A$) if both the map $f:A\to B$ and the map $f^\mathbb N:\mathsf a  \to \mathsf b$ are onto. In this case the image of $f$ is the t-algebra $\mathbf B$.

A \emph{t-isomorphism} is an onto t-homomorphism $f:\mathbf A\to\mathbf B$ such that $f$ is injective. 

If $K$ is a class of t-algebras of the same type, then we denote by $\mathrm I_t K$ the closure of $K$ by t-isomorphism, and 
by $\mathrm H_t K$ the closure of $K$ by t-homomorphic image.

\subsection{t-Congruences}
\begin{definition} Let $\mathbf A$ be a t-algebra of type $\tau$ and trace $\mathsf a$.
   An equivalence relation $\theta$ on $A$ is a \emph{t-congruence} on $\mathbf A$, if  for every $\sigma\in\tau$ and $s,u\in \mathsf a$ we have:
$$s\theta^\mathbb N u\ \text{($s_i\theta u_i$ for every $i$)}\ \Rightarrow \sigma^\mathbf A(s)\theta \sigma^\mathbf A(u).$$
\end{definition}
It is easy to define the quotient t-algebra  $\mathbf A/\theta$.
The kernel of a t-homomorphism is a t-congruence and the map $\pi: A\to A/\theta$ is an onto t-homomorphism of t-algebras.

%
%

\subsection{t-Products} 

Let $\mathbf A_j$ be a  t-algebra of type $\tau$ and trace $\mathsf a_j$ ($j\in J$). Let $B= \Pi_{j\in J} A_j$ and  
 $\mathsf b= \Pi_{j\in J} \mathsf a_j$.
Every $s=(s^j:j\in J)\in \Pi_{j\in J} \mathsf a_j$ is an $\mathbb N\times J$ matrix. If we read $s$ by rows, then $s=(s_k: k\in \mathbb N)\in B^\mathbb N$, where $s_k=(s^j_k:j\in J)\in B$ for every $k\in\mathbb N$.


\begin{definition}
 A t-algebra $\mathbf B$ of type $\tau$ and trace $\mathsf b$ is the \emph{t-product} of the t-algebras $\mathbf A_j$ ($j\in J$)  if $B= \Pi_{j\in J} A_j$, $\mathsf b= \Pi_{j\in J} \mathsf a_j$,
and 
$\sigma^\mathbf B(s)= (\sigma^{\mathbf A_j}(s^j) :j\in J)$ for every $s=(s^j:j\in J)\in \mathsf b$.
\end{definition}

If $K$ is a class of t-algebras of the same type, then we denote by $\mathrm P_t K$ the closure of $K$ by t-product.

\begin{remark}
 The class of all basic t-algebras of type $\tau$ is closed under t-subalgebra, t-homomorphic image and finite t-product. 
 \end{remark}

\section{Functional clone algebras} \label{sec:fca}
In this section we generalise the notion of a functional clone algebra introduced in \cite{BS22}. 
\emph{Functional clone algebras}  are the most natural $\mathsf{CA}$s, the ones the axioms are intended to characterise. 
Given a trace $\mathsf a$ on a given set $A$, the universe of a functional clone algebra is a set of t-operations from $\mathsf a$ into $A$, containing all projections and closed under composition. In this framework  $q_n(\varphi,\psi_1,\dots,\psi_n)$ represents the $n$-ary composition of $\varphi$ with $\psi_1,\dots,\psi_n$ acting on the first $n$ coordinates.

\begin{definition}
 Let $\mathbf A= ( A,\mathsf a,\sigma^\mathbf A)_{\sigma\in \tau}$ be a t-algebra and $A^\mathsf a$ be the set of all t-operations from $\mathsf a$ into $A$.
The tuple $\mathbf A^{\!(\mathsf a)}=(A^\mathsf a, q^\mathsf a_n,\e^\mathsf a_i,\sigma^\mathbf A)_{\sigma\in\tau}$, where, for every $s\in \mathsf a$ and $\varphi,\psi_1,\dots,\psi_n\in  A^\mathsf a$, 
\begin{itemize}
\item $\e_i^\mathsf a(s)=s_i$;
\item $q^\mathsf a_n(\varphi,\psi_1,\dots,\psi_n)(s)=\varphi(s[\psi_1(s),\dots, \psi_n(s)])$,
\end{itemize}
 is called the \emph{full functional clone $\tau$-algebra with value domain $\mathbf A$ and trace $\mathsf a$}.
  A subalgebra of $\mathbf A^{\!(\mathsf a)}$  is called a \emph{functional clone $\tau$-algebra $(\mathsf{FCA}_\tau)$ with value domain $\mathbf A$  and trace $\mathsf a$} ($\mathsf{FCA}_\tau(\mathbf A,\mathsf a)$, for short). 
  \end{definition}
 
   If $\mathsf a$ is a basic (resp. complete) trace, then a subalgebra of $\mathbf A^{\!(\mathsf a)}$  is called a \emph{basic (complete) functional clone $\tau$-algebra with value domain $\mathbf A$ and  trace $\mathsf a$} ($\mathsf{bFCA}_\tau(\mathbf A,\mathsf a)$, for short) (resp. $\mathsf{cFCA}_\tau(\mathbf A,\mathsf a)$).
   
   If $\tau$ is empty, then the full functional pure clone $\tau$-algebra with value domain $A$ and trace $\mathsf a$ is denoted by 
   $A^{(\mathsf a)}=(A^\mathsf a, q^\mathsf a_n,\e^\mathsf a_i)$.

%

\begin{remark}
We warn the reader that in \cite{BS22} the $\mathsf{cFCA}$s were called functional clone algebras, and the $\mathsf{bFCA}$s were called point-relativised functional clone algebras. The $\mathsf{FCA}$s, whose traces are neither basic nor complete,  are introduced for the first time in this article.  
 \end{remark}


The algebraic and functional notions of independence  are equivalent.

\begin{lemma}
 \cite[Lemma 3.8]{BS22} A t-operation  $\varphi:\mathsf a\to A$  is independent of $\e_n$ iff, for all $s,u\in \mathsf a$, $u_i=s_i$ for all $i\neq n$ implies $\varphi(u) = \varphi(s)$.
\end{lemma}

There may exist elements in a clone algebra whose dimension is infinite. For example,  $\mathrm{dim}(\sigma^\mathbf A)=\omega$ for the t-operation $\sigma^\mathbf A$ introduced in Example \ref{exa:diminf}.



The following lemma explains the importance of basic traces.

\begin{lemma}\label{lem:topdim} \cite[Lemma 4.3]{BS22} Let $A$ be a set and  $\mathsf a$ be a basic trace on  $A$. 
  Then the following conditions are equivalent for a t-operation $\varphi:\mathsf a\to A$: 
\begin{enumerate}
\item $\varphi=f_\mathsf a^\top$ is the top $\mathsf a$-extension of a finitary operation $f: A^n\to A$.
\item $\varphi$ has dimension $\leq n$.
\end{enumerate}
Then $\varphi:\mathsf a\to A$ is zero-dimensional iff it is  a constant map.
\end{lemma}

\begin{remark}\label{rem:semi}
 The above equivalence does not hold if the trace is not basic. In general,
a t-operation $\varphi: \mathsf a\to A$ is zero-dimensional iff it is semiconstant (see Definition \ref{def:semi}). Moreover, there exist infinitary operations $\varphi: A^\mathbb N\to A$ of positive finite dimension,  which are not the top extension of any finitary operation (see Section \ref{sec:tafromalg}).
\end{remark}

\subsection{Block algebras}\label{exa:block}   
Let $\mathsf a$ be a trace on a given set $A$ and $B_\mathsf a^\top=\{f_\mathsf a^\top:\mathsf a\to A\ |\ f\in O_A\}$ (see Definition \ref{def:topext}).
 The tuple $\mathcal B_\mathsf a^\top=(B_\mathsf a^\top,q_n^{\mathsf a},\e_i^{\mathsf a})$  is a pure finite-dimensional $\mathsf{FCA}(A,\mathsf a)$.
 
 The $\mathsf{FCA}$ $\mathcal B_\mathsf a^\top$ will be called \emph{the full block algebra on $A$ of trace $\mathsf a$}. A subalgebra of  $\mathcal B_\mathsf a^\top$ is called a \emph{block algebra on $A$}. 
 
 The lattice of all clones
on $A$ is isomorphic to the lattice of all subalgebras of $\mathcal B_\mathsf a^\top$ (see also \cite[Proposition 3.14]{BS22}).
 In \cite[Section 3.3]{BS22}  it was shown that every finite-dimensional clone algebra is isomorphic to a block algebra.

A clone on a given set $A$ determines  a family of t-algebras  and of block algebras on $A$. 

\begin{definition}\label{def:blockalg}
  Let $C$ be a clone of finitary operations on $A$, $\mathsf a$ be a trace on $A$ and $C^\top_\mathsf a=\{f_\mathsf a^\top: f\in C\}$, where $f_\mathsf a^\top: \mathsf a\to A$ is the top $\mathsf a$-extension of $f$. Then $\mathcal C^\top_\mathsf a=(C^\top_\mathsf a, q_n^\mathsf a,\e_i^\mathsf a, f_\mathsf a^\top)_{f\in C}$ is called \emph{the block $C$-algebra of trace $\mathsf a$}, and 
$\mathbf C_A^\top=(A,\mathsf a,f_\mathsf a^\top)_{f\in C}$ is called \emph{the t-algebra of type $C$ and trace $\mathsf a$}.
\end{definition}

\subsection{The t-algebra under a $\mathsf{CA}$} \label{sec:under}

A clone $\tau$-algebra $\mathcal C$ determines a basic trace $\pmb \epsilon  =[\epsilon^\mathcal C]_\mathbb N$, where $\epsilon^\mathcal C = (\e_i^\mathcal C:i\in\mathbb N)$.
For every $c\in C$, we define  a t-operation
$\varphi_c:  \pmb \epsilon \to C$ as follows: 
$$\varphi_c(s)=q_n^\mathcal C(c,s_1,\dots,s_n),\ \text{for every $n$ such that $s=\epsilon^\mathcal C[s_1,\dots,s_n]$.}$$ 
The set $O_\mathcal C=\{\varphi_c: c\in C\}$ is called \emph{the set of all t-operations of $\mathcal C$}.

\begin{definition}\label{def:underC}
  Given a clone $\tau$-algebra $\mathcal C$, the tuple  
  $\mathcal C^\downarrow=( C, \pmb \epsilon,\varphi_{\sigma^\mathcal C})_{\sigma\in\tau}$
is called \emph{the t-algebra  under $\mathcal C$} of type $\tau$ and trace $\pmb \epsilon$.
\end{definition}

\begin{lemma}\label{lem:lemma1} The correspondence $\mathcal C\mapsto \mathcal C^\downarrow$ determines a functor from the category of clone $\tau$-algebras into the category of t-algebras of type $\tau$. 
\end{lemma}

\begin{proof}  Every homomorphism $f:\mathcal C\to \mathcal D$ of clone $\tau$-algebras is also a t-homomorphism from  $\mathcal C^\downarrow$ into $\mathcal D^\downarrow$.
  \[
\begin{array}{lll}
\varphi_{\sigma^\mathcal D}\circ f^\mathbb N(\epsilon^\mathcal C[s_1,\dots,s_n])  &  = & \varphi_{\sigma^\mathcal D}( \epsilon^\mathcal D[f(s_1),\dots,f(s_n)])   \\
  & =  & q_n^\mathcal D(\sigma^\mathcal D,f(s_1),\dots,f(s_n))  \\
  & =  & q_n^\mathcal D(f(\sigma^\mathcal C),f(s_1),\dots,f(s_n))  \\
  &  = &  f(q_n^\mathcal C(\sigma^\mathcal C,s_1,\dots,s_n)) \\
  &  = & f\circ\varphi_{\sigma^\mathcal C}(\epsilon^\mathcal C[s_1,\dots,s_n]). \qedhere
\end{array}
\]
\end{proof}



In the following lemma a clone algebra is shown to be isomorphic to a $\mathsf{bFCA}$.

 \begin{lemma}\label{lem:ddd} \cite[Lemma 4.5]{BS22} Let  $\mathcal C$ be a clone $\tau$-algebra.

(i) The algebra  $\mathcal O_{\mathcal C}=(O_\mathcal C, q_n^{\pmb \epsilon},\e_i^{\pmb \epsilon}, \varphi_{\sigma^\mathcal C})_{\sigma\in\tau}$ is a $\mathsf{bFCA}_{\tau}(\mathcal C^\downarrow, \pmb \epsilon)$.

 (ii) $\mathcal C\cong \mathcal O_{\mathcal C}$.
 \end{lemma}

\begin{proof}
The map $c\in C \mapsto \varphi_c\in O_\mathcal C$ is an isomorphism from $\mathcal C$ onto $\mathcal O_{\mathcal C}$.
\end{proof}

$\mathcal O_{\mathcal C}$ is  called \emph{the clone algebra of all t-operations of  $\mathcal C$}.

\begin{remark}\label{rem:neu-clo}
 Definition \ref{def:underC} and Lemma \ref{lem:ddd} justify the notions of trace and t-algebra.
 To develop universal clone algebra we have two possibilities: we can use either  an infinitary approach based on Neumann's abstract $\aleph_0$-clones or a finitary approach based on clone algebras. In the first case, we do not need traces, because,  for every abstract $\aleph_0$-clone $\mathcal C=(C,q_\infty^\mathcal C,\e_i^\mathcal C)$ and every $c\in C$, we can define  an infinitary operation
$\phi_c:  C^\mathbb N \to C$ as follows: 
$$\phi_c(s)=q_\infty^\mathcal C(c,s_1,\dots,s_n,\dots),\ \text{for every $s\in C^\mathbb N$.}$$ 
The map $c\mapsto \phi_c$ determines an isomorphism from $\mathbf C$ onto a functional $\aleph_0$-clone with value domain $C$ (see Definition \ref{sec:neu}).
The price to be paid is that the abstract $\aleph_0$-clones are not finitary algebras. 
We remark that in one of the main results of \cite{BS22} it was shown that every clone $\tau$-algebra $\mathcal C$  is isomorphic to a complete $\mathsf{FCA}_\tau$, whose value domain is related in a non-trivial way to $C$.
\end{remark}

\section{Term clone algebras and hyperidentities}\label{sec:cats}
In this section we introduce the clone $\tau$-algebra $\mathcal T_\tau(X)$ of $\tau$-hyperterms and show that it is the free algebra over the set $X$ of generators in the variety of all clone $\tau$-algebras.  The $\tau$-hyperterms  over $\emptyset$ are called $\tau$-terms.
The t-algebra $\mathcal T_\tau(\emptyset)^\downarrow$ under $\mathcal T_\tau(\emptyset)$ will be shown to be the absolutely free t-algebra over $\epsilon=(\e_1,\dots,\e_n,\dots)$ in the class of all t-algebras.
We also introduce the term clone $\tau$-algebra $\mathbf A^\uparrow$  over a t-algebra $\mathbf A$. It is the analogue  of 
the term clone of an algebra.

\begin{definition}\label{def:hyper} Let $\tau$ be a set of operator symbols and $X$ be an infinite set such that $X\cap \tau=\emptyset$.
  The set $T_{\tau}(X)$ of  \emph{$\tau$-hyperterms over $X$}  is built up by induction as follows:
\begin{enumerate}
\item $\e_1,\dots,\e_n,\dots$ are $\tau$-hyperterms;
\item If $t_1,\dots,t_n$ are $\tau$-hyperterms 
and $w\in \tau\cup X$, then $w(t_1,\dots,t_n,\e_{n+1},\e_{n+2},\dots)$ is a $\tau$-hyperterm, for every $n\geq 0$.
\end{enumerate}
\end{definition}
The $\tau$-hyperterms  over $\emptyset$ will be called \emph{$\tau$-terms}. $T_\tau$ denotes the set of all $\tau$-terms.

If there is no danger of ambiguity, for every $w\in\tau\cup X$, we write  $w$  for $w(\e_1,\dots,\e_n,,\dots)$.
We identify $X(\bar\e)=\{w(\e_1,\e_2,\dots,\e_n,\dots) : w\in X\}$ with $X$ and $\tau(\bar\e)=\{\sigma(\e_1,\e_2,\dots,\e_n,\dots) : \sigma\in \tau\}$ with $\tau$.


 

\begin{proposition}\label{lem:term}
The algebra $\mathcal{T}_{\tau}(X)=(T_{\tau}(X),q_n^\mathcal T,\e_i^\mathcal T,\sigma^\mathcal T)_{\sigma\in\tau}$, where
 $\e_i^\mathcal T=\e_i$;  $\quad\sigma^\mathcal T=\sigma$ for every $\sigma\in\tau$; and $q_n^\mathcal T$ ($n\geq 0$) is defined by induction as follows:
\begin{enumerate}
\item[(i)] $q_n^\mathcal T(\e_i,t_1,\dots,t_n)=t_i$ $(1\leq i\leq n)$;
\item[(ii)] $q_n^\mathcal T(\e_i,t_1,\dots,t_n)=\e_i$ $(i>n)$;
\item[(iii)] $q_n^\mathcal T(w(t_1,\dots,t_k,\e_{k+1},\e_{k+2},\dots),\bar u)= $\\
$ w(q_n^\mathcal T(t_1,\bar u),\dots, q_n^\mathcal T(t_k,\bar u),q_n^\mathcal T(\e_{k+1},\bar u),q_n^\mathcal T(\e_{k+2},\bar u),\dots)$ for every $w\in\tau\cup X$;
\end{enumerate}
 is the free clone $\tau$-algebra over  the set $X$ of generators in the variety $\mathsf{CA}_\tau$ of clone $\tau$-algebras. If $X=\emptyset$, then $\mathcal{T}_{\tau}$ is initial in $\mathsf{CA}_\tau$.
\end{proposition}

\begin{proof} Let  $\bar\e=\e_1,\dots,\e_n$, $\bar t=t_1,\dots,t_k,\dots$, where $t_k=\e_k$ for almost all $k$, and $w\in \sigma\cup X$.
  \begin{enumerate}
\item[(C3)]   
$q_n^\mathcal T(w(\bar t),\bar\e)=w(q_n^\mathcal T(t_1,\bar \e),\dots, q_n^\mathcal T(t_k,\bar \e),\dots)= w(\bar t)$.
\item[(C4)] 
$
q_{n+r}^\mathcal T(w(\bar t),\mathbf y,\e_{n+1},\dots,\e_{n+r})
=w(q_{n+r}^\mathcal T(t_1,\mathbf y,\e_{n+1},\dots,\e_{n+r}),\dots)=$\\ $w(q_n^\mathcal T(t_1,\mathbf y),\dots)
=q_n^\mathcal T(w(\bar t),\mathbf y).$\\
\item[(C5)] 
 $
\begin{array}{lll}
&& q_n^\mathcal T(q_n^\mathcal T(w(\bar t),\mathbf y), \mathbf z)\\
   &  = &q_n^\mathcal T(w(q_n^\mathcal T(t_1,\mathbf y),\dots, q_n^\mathcal T(t_k,\mathbf y),\dots),\mathbf z)\\
&  = &w(q_n^\mathcal T(q_n^\mathcal T(t_1,\mathbf y),\mathbf z),\dots,q_n^\mathcal T(q_n^\mathcal T(t_k,\mathbf y),\mathbf z),\dots)\\
&  = &w(q_n^\mathcal T(t_1, q_n^\mathcal T(y_1,\mathbf z),\dots, q_n^\mathcal T(y_n,\mathbf z)),\dots,\\
&& \qquad \qquad \qquad\qquad q_n^\mathcal T(t_k, q_n^\mathcal T(y_1,\mathbf z),\dots, q_n^\mathcal T(y_n,\mathbf z)),\dots)\\
   &  = &  q_n^\mathcal T(w(\bar t),q_n^\mathcal T(y_1,\mathbf z),\dots,q_n^\mathcal T(y_n,\mathbf z)).
 \end{array}
$
\end{enumerate}
  Let $\mathcal C$ be a clone $\tau$-algebra and $\alpha:\tau\cup X\to C$ be a function such that $\alpha(\sigma)=\sigma^\mathcal C$. We prove that there exists a unique homomorphism $\alpha^\star :\mathcal T_{\tau}(X) \to \mathcal C$ such that $\alpha^\star(w)=\alpha(w)$ for every $w\in X$: 
 $$\alpha^\star(\e_i^\mathcal T)=\e_i^\mathcal C;\quad
 \alpha^\star(w(t_1,\dots,t_k,\e_{k+1},\e_{k+2},\dots))= q_k^\mathcal C(\alpha(w),\alpha^\star(t_1),\dots,\alpha^\star(t_k)).$$
We prove that $\alpha^\star$ is a homomorphism. We write $w(\bar t)$ for $w(t_1,\dots,t_k,t_{k+1},\dots)$ and $\alpha^\star(\bar u)$ for $\alpha^\star(u_1),\dots,\alpha^\star(u_n)$.
\[
\begin{array}{llll}
  &   & \alpha^\star(q_n^\mathcal T(w(\bar t),\bar u)) & \\
  &=   &\alpha^\star(w(q_n^\mathcal T(t_1,\bar u),\dots, q_n^\mathcal T(t_k,\bar u),q_n^\mathcal T(t_{k+1},\bar u),\dots)) &  \\
  & =  &   q_m^\mathcal C(\alpha(w),\alpha^\star(q_n^\mathcal T(t_1,\bar u)),\dots,\alpha^\star(q_n^\mathcal T(t_m,\bar u)))\quad\text{($m>k,n$)}&\\
  & =  &  q_m^\mathcal C(\alpha(w),q_n^\mathcal C(\alpha^\star(t_1),\alpha^\star(\bar u)),\dots,q_n^\mathcal C(\alpha^\star(t_m),\alpha^\star(\bar u))) &\\
& =  & q_m^\mathcal C(\alpha(w),q_m^\mathcal C(\alpha^\star(t_1),\alpha^\star(\bar u),\e_{n+1},\dots,\e_m),\dots,&\\
&&\qquad\qquad\qquad\qquad\qquad\qquad q_m^\mathcal C(\alpha^\star(t_m),\alpha^\star(\bar u),\e_{n+1},\dots,\e_m)) &\\
& =  &q_m^\mathcal C(q_m^\mathcal C(\alpha(w),\alpha^\star(t_1),\dots,\alpha^\star(t_m)),\alpha^\star(\bar u),\e_{n+1},\dots,\e_m)& \\
& =  &q_n^\mathcal C(q_m^\mathcal C(\alpha(w),\alpha^\star(t_1),\dots,\alpha^\star(t_m)),\alpha^\star(\bar u))& \\
& =  &q_n^\mathcal C(\alpha^\star(w(\bar t)),\alpha^\star(\bar u)). 
\end{array}
\]
Let $g:\mathcal T_{\tau}(X) \to \mathcal C$ be another homomorphism such that $g(w)=\alpha(w)$. Then we have:\\
$g(w(t_1,\dots,t_k,\e_{k+1},\dots))= g(q_k^\mathcal T(w,t_1,\dots,t_k))=
q_k^\mathcal C(g(w),g(t_1),\dots,g(t_k) ) =q_k^\mathcal C(\alpha(w),\alpha^\star(t_1),\dots,\alpha^\star(t_k) )=\dots=\alpha^\star(w(t_1,\dots,t_k,\e_{k+1},\dots))$.\qedhere
\end{proof}

The t-algebra $\mathcal T_\tau^\downarrow$ under $\mathcal T_\tau$ will be called the \emph{t-algebra of $\tau$-terms}.   
 The t-algebra $\mathcal T_\tau(X)^\downarrow$ under $\mathcal T_\tau(X)$ will be called the \emph{t-algebra of $\tau$-hyperterms}.

\subsection{The term clone algebra over a t-algebra} \label{sec:over}

In this section we introduce the term clone algebra $\mathbf A^\uparrow$ and the polynomial clone algebra $\mathbf A^{\!\Uparrow}$ over a t-algebra $\mathbf A$.   We define three equational theories of $\mathbf A$: the equational t-theory, the equational w-hypertheory and the equational f-hypertheory.
We prove some propositions and lemmas will be useful in the following.


\begin{definition}
 Let $\mathbf A$ be a t-algebra of type $\tau$ and trace $\mathsf a$. 
 \begin{itemize}
\item  The \emph{term clone $\tau$-algebra $\mathbf A^\uparrow$  over $\mathbf A$} is the minimal subalgebra of the full $\mathsf{FCA}_\tau$  $\mathbf A^{\!(\mathsf a)}$.
\item The \emph{polynomial clone $\tau$-algebra $\mathbf A^{\!\Uparrow}$  over $\mathbf A$} is the minimal subalgebra of the full $\mathsf{FCA}_\tau$  $\mathbf A^{\!(\mathsf a)}$ containing all constant t-operations. 
\end{itemize}
\end{definition}

If $K$ is a class of t-algebras, then we denote by $K^{\!\uparrow}$ the class $\{\mathbf A^{\!\uparrow} : \mathbf A\in K\}$, and by $K^{\!\Uparrow}$ the class $\{\mathbf A^{\!\Uparrow} : \mathbf A\in K\}$.

$\mathbf A^\uparrow$ is the image of the unique  homomorphism $(-)^\mathbf A:\mathcal T_\tau\to \mathbf A^{\!(\mathsf a)}$ mapping $t\in T_\tau$ into  $t^\mathbf A:\mathsf a\to A$. The term operation $t^\mathbf A$ is defined by induction as follows, for every $s\in \mathsf a$:
$$t^\mathbf A(s)=\begin{cases}
s_i&\text{if $t\equiv \e_i$}\\
\sigma^\mathbf A(t_1^\mathbf A(s),\dots,t_n^\mathbf A(s),s_{n+1},\dots)&\text{if $t\equiv \sigma(t_1,\dots,t_n,\e_{n+1},\dots)$}
\end{cases}$$

We remark that the t-algebra $\mathbf A^{\updownarrows}$ under the term clone algebra $\mathbf A^\uparrow$ is not in general isomorphic to $\mathbf A$. In Section \ref{sec:free} we introduce some conditions that imply $\mathbf A^{\updownarrows}\cong \mathbf A$ (see also Proposition \ref{prop:thm}).

\begin{example}\label{lem:coincide}
 Let $C$ be a clone on $A$ and $\mathbf C_\mathsf a^\top$ be the t-algebra of type $C$ and trace $\mathsf a$ (see Definition \ref{def:blockalg}). Then the term clone algebra  $(\mathbf C_\mathsf a^\top)^\uparrow$ over  $\mathbf C_\mathsf a^\top$ is equal to the block $C$-algebra $\mathcal C^\top_\mathsf a$ of trace $\mathsf a$ introduced in Definition \ref{def:blockalg}.
\end{example}

\subsection{Identities and hyperidentities} 

A $\tau$-identity is a pair $(t_1,t_2)$, written $t_1=t_2$,  of $\tau$-terms. A $\tau$-hyperidentity is a pair $(u_1,u_2)$, written $u_1=u_2$,  of $\tau$-hyperterms. 

Recall from Definition \ref{def:semi} the notion of semiconstant t-operation.

A homomorphism $\alpha:\mathcal T_\tau(X)\to\mathbf A^{\!(\mathsf a)}$  is \emph{weak} (resp. \emph{strong}) 
if  $\alpha(x)$ is a constant  (resp. semiconstant) t-operation, 
for every $x\in X$. In the weak case, the image of $\alpha$ is a subalgebra of $\mathbf A^{\!\Uparrow}$. 

\begin{definition}\label{def:sat} Let  Let  $t_1,t_2\in T_\tau$ be $\tau$-terms, $u_1,u_2\in T_\tau(X)$ be $\tau$-hyperterms, and $\mathbf A$ be a t-algebra of type $\tau$ and trace $\mathsf a$. We say that 
\begin{enumerate}
\item \emph{$\mathbf A$ satisfies $t_1=t_2$}, and we write $\mathbf A\models t_1=t_2$, if $t_1^\mathbf A=t_2^\mathbf A$. The kernel $\theta_\mathbf A$ of  the unique  homomorphism $(-)^\mathbf A:\mathcal T_\tau\to \mathbf A^{\!\uparrow}$ is called the \emph{equational t-theory of $\mathbf A$}.

\item  $\mathbf A$ \emph{weakly satisfies} the $\tau$-hyperidentity $u_1=u_2$, and we write $\mathbf A\models_w u_1=u_2$, 
if $\alpha(u_1)=\alpha(u_2)$ for every weak homomorphism $\alpha:\mathcal T_\tau(X)\to\mathbf A^{\!(\mathsf a)}$.
$\mathrm{Tw}^X_{\mathbf A}=\{u_1=u_2 : \mathbf A\models_w u_1=u_2\}$ is called the \emph{equational w-hypertheory of $\mathbf A$}.

\item $\mathbf A$ \emph{fully satisfies} the $\tau$-hyperidentity $u_1=u_2$, and we write $\mathbf A\models_f u_1=u_2$, if $\alpha(u_1)=\alpha(u_2)$ for every homomorphism $\alpha:\mathcal T_\tau(X)\to\mathbf A^{\!(\mathsf a)}$. $\mathrm{Tf}^X_{\mathbf A}=\{u_1=u_2 : \mathbf A\models_f u_1=u_2\}$ is called the \emph{equational f-hypertheory of $\mathbf A$}.

\end{enumerate}
\end{definition}

If $K$ is a class of t-algebras, we define:
\begin{itemize}
\item $\theta_K = \bigcap_{\mathbf A\in K} \theta_\mathbf A$ is the equational t-theory of $K$.
\item $\mathrm{Tw}^X_{K} = \bigcap_{\mathbf A\in K}\mathrm{Tw}^X_{\mathbf A}$ is called the equational w-hypertheory of $K$.
\item $\mathrm{Tf}^X_{K} = \bigcap_{\mathbf A\in K} \mathrm{Tf}^X_{\mathbf A}$ is called the equational f-hypertheory of $K$.
\end{itemize}
If there is no danger of confusion, we omit $X$ in $\mathrm{Tf}^X_{K}$ and $\mathrm{Tw}^X_{K}$.

\bigskip

A t-algebra $\mathbf A$ \emph{strongly satisfies} the $\tau$-hyperidentity $u_1=u_2$, and we write $\mathbf A\models_s u_1=u_2$, 
if $\alpha(u_1)=\alpha(u_2)$ for every strong homomorphism $\alpha:\mathcal T_\tau(X)\to\mathbf A^{\!(\mathsf a)}$.

\begin{proposition}\label{prop:ws} $\mathbf A\models_w u_1=u_2$ iff $\mathbf A\models_s u_1=u_2$.
\end{proposition}

\begin{proof} ($\Rightarrow$) 
 Let $\mathbf A\models_w u_1=u_2$,  and $\beta:\mathcal T_\tau(X)\to\mathbf A^{\!(\mathsf a)}$ be a strong homomorphism. We show that $\beta(u_1)=\beta(u_2)$. For every $\mathsf b\subseteq_{\mathrm{bas}} \mathsf a$ and $x\in X$, $\beta(x)_{|\mathsf b}$ is a constant function. Let $\gamma_\mathsf b:X\to A$ be defined as follows: $\gamma_\mathsf b(x) = \beta(x)(s)$ for some (and then all) $s\in\mathsf b$.
Define $\alpha_\mathsf b: \mathcal T_\tau(X)\to\mathbf A^{\!(\mathsf a)}$ be the unique weak homomorphism such that $\alpha_\mathsf b(x)(s) =\gamma_\mathsf b(x)$ for every $s\in\mathsf a$. Then $\alpha_\mathsf b(u_1)=\alpha_\mathsf b(u_2)$, for every $\mathsf b\subseteq_{\mathrm{bas}} \mathsf a$.
It is easy to prove by induction on $w\in T_\tau(X)$ that $\beta(w)(s)=\alpha_\mathsf b(w)(s)$, for every $s\in \mathsf b$ and $\mathsf b\subseteq_{\mathrm{bas}} \mathsf a$. Then $\beta(u_1)(s)= \alpha_\mathsf b(u_1)(s)=\alpha_\mathsf b(u_2)(s)= \beta(u_1)(s)$ for every $s\in\mathsf b$ and every $\mathsf b\subseteq_{\mathrm{bas}} \mathsf a$.
\end{proof}

\subsection{Functional clone algebras and t-algebras}
The following proposition describes the existing relations between a t-algebra $\mathbf A$ and the functional clone algebras of value domain $\mathbf A$. 

\begin{proposition}\label{prop:zzz} 
Let $\mathbf A$ be a t-algebra of type $\tau$ and trace $\mathsf a$,  $\mathcal F$ be a $\mathsf{FCA}_\tau(\mathbf A,\mathsf a)$ and $s\in\mathsf a$. Then the following conditions hold.
 
  (i) The function $\bar s: F\to A$, defined by $\bar s(\varphi)=\varphi(s)$, is a t-homomorphism from $\mathcal F^\downarrow$ into $\mathbf A$.

(ii) $\mathcal F^\downarrow\in \mathrm S_t\mathrm P_t(\mathbf A)$. 

(iii) The function $\hat{s}:T_\tau\to A$ mapping $t\mapsto t^\mathbf A(s)$, is a t-homomorphism from $\mathcal T_\tau^\downarrow$ into $\mathbf A$, whose kernel $\theta_{\hat s}$ satisfies
 $\theta_\mathbf A\subseteq \theta_{\hat s}$.
 
 (iv)  $\theta_\mathbf A= \bigcap_{r\in \mathsf a} \theta_{\hat r}$. 

\end{proposition}

\begin{proof} (i) We recall that $\mathcal F^\downarrow$ is a t-algebra of trace $\pmb \epsilon=[\epsilon^\mathcal F]_\mathbb N$. The map $\bar s^\mathbb N:  \pmb \epsilon  \to \mathsf a$ is defined as follows for every $u=\epsilon^\mathcal F[\psi_1,\dots,\psi_k]\in   \pmb \epsilon$: 
$\bar s^\mathbb N(u)=s[\psi_1(s),\dots,\psi_k(s)]\in \mathsf a$.
Since $\varphi_{\sigma^\mathbf A}(u)=_{\mathrm{def}} q_k^{\mathsf a}(\sigma^\mathbf A,\psi_1,\dots,\psi_k)$, then
$$\bar s(\varphi_{\sigma^\mathbf A}(u))=\varphi_{\sigma^\mathbf A}(u)(s)= \sigma^\mathbf A(s[\psi_1(s),\dots,\psi_k(s)])= \sigma^\mathbf A(\bar s^\mathbb N(u)).$$

(ii) The product $\langle \bar s: s\in\mathsf a\rangle$ of the t-homomorphisms $\bar s$ embeds  $\mathcal F^\downarrow$ into the t-power $\mathbf A^\mathsf a$ ($\mathsf a$ copies of $\mathbf A$).

 (iii) By Lemma \ref{lem:lemma1}  $($-$)^\mathbf A: \mathcal T_\tau^\downarrow \to \mathbf A^\updownarrows$ is a t-homomorphism. Similarly,  by (i) $\bar s: \mathbf A^\updownarrows\to\mathbf A$ is also a t-homomorphism. Then we define $\hat s=\bar s\circ ($-$)^\mathbf A$.
%
\end{proof}


\begin{proposition}\label{prop:natural}
 If $f:\mathbf A\to \mathbf B$ is an onto t-homomorphism of t-algebras, then the correspondence $t^\mathbf A\mapsto t^\mathbf B$  is a homomorphism from $\mathbf A^\uparrow$ onto $\mathbf B^\uparrow$.
\end{proposition}

\begin{proof} Let $\mathsf a$ (resp. $\mathsf b$) be the trace of $\mathbf A$ (resp. $\mathbf B$).
It is sufficient to prove that, if $t_1^\mathbf A=t_2^\mathbf A$, then $t_1^\mathbf B=t_2^\mathbf B$.  By Proposition \ref{prop:zzz}
 $f\circ \bar s: \mathbf A^{\updownarrows}\to \mathbf B$, where $s\in \mathsf a$, is a t-homomorphism satisfying $f(t^\mathbf A(s))=t^\mathbf B(f^\mathbb N (s))$. 
Since $f^\mathbb N: \mathsf a\to \mathsf b$ is onto, it follows that if $t_1^\mathbf A=t_2^\mathbf A$, then $t_1^\mathbf B=t_2^\mathbf B$.
\end{proof}

\subsection{Minimal clone algebras}

\begin{definition}\label{def:min}
 A clone $\tau$-algebra $\mathcal C$ is \emph{minimal} if it is generated by the constants $\e_i^\mathcal C$ ($i \geq 1$) and $\sigma^\mathcal C$ ($\sigma\in\tau$).
\end{definition}
 
\begin{proposition}\label{prop:thm} Let $\mathcal D$ be a  clone $\tau$-algebra.
\begin{enumerate}
\item If $\mathcal D$ is minimal, then $\mathcal D^{\downuparrows} \cong \mathcal D$. 
\item $(\mathcal T_\tau/\phi)^{\downuparrows} \cong \mathcal T_\tau/\phi$, for every congruence $\phi$ of $ \mathcal T_\tau$.
\item $\mathbf A^{\uparrow\downuparrows} \cong \mathbf A^\uparrow$, for every t-algebra $\mathbf A$.
\end{enumerate}
\end{proposition}

\begin{proof}
(i) By Lemma \ref{lem:ddd}  $\mathcal D$ is isomorphic to the $\mathsf{FCA}$ $\mathcal O_{\mathcal D}$ with value domain  $\mathcal D^\downarrow$. Since $\mathcal D$ is minimal then  $ \mathcal O_{\mathcal D}$ is also minimal; hence, $\mathcal O_{\mathcal D}$ is $\mathcal D^{\downuparrows}$. 
(ii) and (iii) are trivial consequences of (i).
\end{proof}

%
%
%
%

\section{t-Varieties, Ft-varieties and Et-varieties}\label{sec:tvarieties}
In this section we define t-varieties, Ft-varieties and Et-varieties of t-algebras. A class of t-algebras is a t-variety if it is closed under $\mathrm H_t$, $\mathrm S_t$ and $\mathrm P_t$. An Et-variety (Ft-variety) is a t-variety closed under (full) t-expansion. 

In the following $\mathbf A$ is a t-algebra of type $\tau$ and trace $\mathsf a$.

\subsection{Sum of t-algebras and (F)t-varieties}

Let $\mathsf a$ be a trace on a given set $A$. A  \emph{t-partition of $\mathsf a$} is a family $\{\mathsf b_i\}_{i\in I}$ of traces on $A$ such that
$\bigcup_{i\in I}\mathsf b_i = \mathsf a$ and $\mathsf b_i\cap \mathsf b_j = \emptyset$ for every $i\neq j$.
The family $\{\mathsf b_i\}_{i\in I}$ is a basic t-partition if $\mathsf b_i$ is a basic trace for every $i\in I$.


Recall from Definition \ref{def:fullsub} the notion of full t-subalgebra.

\begin{definition}
 Let $\mathbf A=(A,\mathsf a,\sigma^\mathbf A)_{\sigma\in\tau}$ be a t-algebra. We say that
$\mathbf A$ is \emph{the sum of the full t-subalgebras $\mathbf A_{\upharpoonright\mathsf b_i}$ ($i\in I$) of $\mathbf A$}, and we write $\mathbf A=\oplus_{i\in I} \mathbf A_{\upharpoonright\mathsf b_i}$,  if $\{\mathsf b_i\}_{i\in I}$ is a t-partition of $\mathsf a$.
\end{definition}


\noindent There exists a unique basic t-partition of $\mathsf a$. We write $\oplus^b_{i\in I} \mathbf A_{\upharpoonright\mathsf b_i}$ for $\oplus_{i\in I} \mathbf A_{\upharpoonright\mathsf b_i}$.

\begin{definition}
  We say that a class $K$ of t-algebras of type $\tau$ \emph{is closed under full t-expansion} (and we write $\mathrm F_t K= K$) if, for every t-algebra $\mathbf A=\oplus^b_{i\in I} \mathbf A_{\upharpoonright\mathsf b_i}$ of type $\tau$ and trace $\mathsf a$,
$$(\forall i\in I.\ \mathbf A_{\upharpoonright\mathsf b_i}\in K) \Rightarrow \mathbf A\in K.$$
\end{definition}

\begin{definition} Let $K$ be a class of t-algebras of type $\tau$.

(i) $K$ is a \emph{t-variety} if it is closed under $\mathrm H_t$, $\mathrm S_t$ and $\mathrm P_t$.

(ii) $K$ is an \emph{Ft-variety} if it is a t-variety closed under $F_t$.
\end{definition}

\begin{proposition}\label{prop:c(aa)} Let $\mathbf A=\oplus_{i\in I} \mathbf A_{\upharpoonright\mathsf b_i}$.

(i) There exists a unique homomorphism from $\mathbf A^\uparrow$ onto $\mathbf A_{\upharpoonright\mathsf b_i}^\uparrow$, defined by
 $t^\mathbf A\mapsto t^{\mathbf A}_{|\mathsf b_i}$.
 
(ii) $\mathbf A^\uparrow$ is a subdirect product of the term clone algebras $\mathbf A_{\upharpoonright\mathsf b_i}^\uparrow$ ($i\in I$).

(iii) $\theta_\mathbf A= \bigcap_{i\in I} \theta_{\mathbf A_{\upharpoonright\mathsf b_i}}$.
\end{proposition}

\subsection{Et-varieties}

\begin{definition}\label{def:As}
 \emph{The t-subalgebra $\mathbf A_{\bar s}$ of $\mathbf A$ generated by $s\in \mathsf a$} is  defined as the image of the t-homomorphism (defined in Proposition \ref{prop:zzz}(i)) $\bar s:\mathbf A^\updownarrows\to\mathbf A$
mapping $t^\mathbf A$ into $t^\mathbf A(s)$. The universe  of $\mathbf A_{\bar s}$ is the set $A_{\bar s} = \{t^\mathbf A(s): t\in T_\tau\}$ and its trace $\mathsf a_{\bar s}= [s]_\mathbb N^A \cap (A_{\bar s})^\mathbb N$ is a basic trace.
\end{definition}

The elements of $A_{\bar s}$ are defined by induction as follows:
\begin{enumerate}
\item $s_i\in A_{\bar s}$ for every $i\in \mathbb N$;
\item If $a_1,\dots,a_n\in A_{\bar s}$ then $\sigma^\mathbf A(a_1,\dots,a_n,s_{n+1},\dots)\in A_{\bar s}$ for every $\sigma\in\tau$ and $n\geq 0$.
\end{enumerate}


The t-subalgebra $\mathbf A_{\bar s}$ is  finitely generated if
$|\mathrm{set}(s)|<\omega$, otherwise, it has a countable infinite set of generators.

%



\begin{proposition}\label{prop:c(a)}
(i) There exists a unique homomorphism from $\mathbf A^\uparrow$ onto $\mathbf A_{\bar s}^\uparrow$, defined by
 $t^\mathbf A\mapsto t^{\mathbf A}_{|\mathsf a_{\bar s}}$.
 
(ii) $\mathbf A^\uparrow$ is a subdirect product of the term clone algebras $\mathbf A_{\bar s}^\uparrow$ ($s\in \mathsf a$).

(iii) $\theta_\mathbf A= \bigcap_{s\in \mathsf a} \theta_{\mathbf A_{\bar s}}$.
\end{proposition}

\begin{proof} (i) Since  $\bar s:\mathbf A^\updownarrows\to\mathbf A_{\bar s}$ is an onto t-homomorphism and by Proposition \ref{prop:thm}(iii)  $\mathbf A^{\updownarrows\uparrow}=\mathbf A^\uparrow$, then by Proposition \ref{prop:natural}
the map $t^\mathbf A \mapsto t^{\mathbf A_{\bar s}}= t^{\mathbf A}_{|\mathsf a_{\bar s}}$ is a homomorphism of clone algebras.
 The other items are a trivial consequence of (i).
\end{proof}

\begin{definition}
  We say that a class $K$ of t-algebras of type $\tau$ \emph{is closed under expansion} (and we write $\mathrm E_t K= K$) if, for every t-algebra $\mathbf A$ of type $\tau$ and trace $\mathsf a$,
$$(\forall s\in \mathsf a.\ \mathbf A_{\bar s}\in K) \Rightarrow \mathbf A\in K.$$
\end{definition}

\begin{definition} Let $K$ be a class of t-algebras of type $\tau$.
$K$ is an  \emph{Et-variety} if it is a t-variety closed under $E_t$.

\end{definition}



\section{Up-and-down between classes of clone algebras and t-algebras} \label{sec:updown}
We have two algebraic levels. 
The lower degree of t-algebras and the higher degree of clone algebras. There are many ways to collectively move between these levels.
 If $K$ is a class of t-algebras, then $K^\smalltriangleup$ is a class of clone algebras. If $H$ is a class of clone algebras, we have two ways to go down: $H^\smalltriangledown$ and  $H^\blacktriangledown$.
In this section we prove that
\begin{itemize}
\item  If $K$ is a t-variety of t-algebras, then $K^\smalltriangleup$ is a  variety of clone algebras.
\item  If $H$ is a class of clone algebras, then  $H^\smalltriangledown$ is an Ft-variety and  $H^\blacktriangledown$ is an Et-variety of t-algebras.
\end{itemize}

 A $\mathsf{FCA}_\tau(\mathbf A,\mathsf a)$ $\mathcal F$ is \emph{rich} if the polynomial clone $\tau$-algebra $\mathbf A^{\!\Uparrow}$ is a subalgebra of $\mathcal F$.
$\mathsf{rFCA}$ (resp. $\mathsf{rFCA}(\mathbf A)$) denotes the class of all rich $\mathsf{FCA}$s (with value domain $\mathbf A$).

In this section $K$ is a class of t-algebras of type $\tau$ and $H$ is a class of clone $\tau$-algebras.
We define 

\begin{itemize}
\item $K^\smalltriangleup= \mathbb I \{\mathcal F :\ \text{$\mathcal F$ is a $\mathsf{FCA}_\tau$ with value domain $\mathbf A\in K$}\}$

\item $H^\smalltriangledown=\{\mathbf A: \mathsf{rFCA}_\tau(\mathbf A)\cap  H\neq\emptyset  \}$. 
\item $H^\blacktriangledown=\{\mathbf A: \mathsf{FCA}_\tau(\mathbf A)\cap  H\neq\emptyset  \}$.
\end{itemize}

We have $H^\smalltriangledown\subseteq  H^\blacktriangledown$ and $K\subseteq K^{\smalltriangleup\smalltriangledown}\subseteq K^{\smalltriangleup\blacktriangledown}$.



\begin{lemma}\label{sss} Let $\mathcal F_i$ be a $\mathsf{FCA}_\tau(\mathbf A_i,\mathsf a_i)$  ($i\in I$).
 Then the product $\Pi_{i\in I} \mathcal F_i$ can be embedded into
the full  $\mathsf{FCA}_{\tau}$ with value domain $\Pi_{i\in I}  \mathbf A_i$. If $\mathcal F_i$ is rich for every $i\in I$, then the image of the product $\Pi_{i\in I} \mathcal F_i$ is also rich.
\end{lemma}

\begin{proof} Let $\mathsf b=\Pi_{i\in I} \mathsf a_i$.
The  sequence $\langle \varphi_i: \mathsf a_i \to A_i\in F_i\ |\ i\in I\rangle \in \Pi_{i\in I}  F_i$ maps to $\varphi : \mathsf b \to \Pi_{i\in I}  A_i$, where, for every $s=(s^i:i\in I)\in  \mathsf b$, $\varphi(s)= \langle \varphi_i(s^i)\ |\ i\in I\rangle$. 
\end{proof}

%

\begin{theorem}\label{thm:ggg}
  If $K$ is a t-variety, then $K^\smalltriangleup$  is a variety of clone $\tau$-algebras.
\end{theorem}

\begin{proof} By definition $K^\smalltriangleup$ is  closed under subalgebra and
by Lemma \ref{sss} it is  closed under product. 
 Let $f: \mathcal F\to \mathcal G$ be an onto homomorphism of clone $\tau$-algebras, where   $\mathcal F\in K^\smalltriangleup$ is a $\mathsf{FCA}_{\tau}(\mathbf A)$ with $\mathbf A\in K$. We have to show that $\mathcal G\in K^\smalltriangleup$.
   By Proposition \ref{prop:zzz}(ii)  $\mathcal F^\downarrow \in\mathrm S_t\mathrm P_t(\mathbf A)\subseteq K$.  By Lemma \ref{lem:lemma1} $f$ is also an onto t-homomorphism from $\mathcal F^\downarrow$ onto $\mathcal G^\downarrow$. Since  $\mathcal F^\downarrow\in K$ and $K$ is a t-variety, then  $\mathcal G^\downarrow\in K$.
By Lemma \ref{lem:ddd} $\mathcal G$ is isomorphic to a $\mathsf{FCA}_{\tau}(\mathcal G^\downarrow)$.  It follows that $\mathcal G\in K^\smalltriangleup$.
\end{proof}

\begin{lemma}\label{lem:11} Let $\mathbf A=(A,\mathsf a,\sigma^\mathbf A)$ and $\mathbf B=(B,\mathsf b,\sigma^\mathbf B)$ be t-algebras of type $\tau$,   $f: \mathbf A\to \mathbf B$ be an onto t-homomorphism  and
  $F=\{\varphi: \mathsf a \to A\ |\ \exists \psi:\mathsf b\to B.\ \psi\circ f^\mathbb N=f\circ \varphi\}$. Then  the following conditions hold:
\begin{enumerate}
\item $F$  is a rich subalgebra of the full $\mathsf{FCA}$ $\mathbf A^{\!(\mathsf a)}$.
\item The map $f^\bullet: F\to B^\mathsf b$, defined by 
$f^\bullet(\varphi)(f^\mathbb N(s))=f(\varphi(s))$ for every $s\in \mathsf a$ and $\varphi\in F$,
is a homomorphism from the clone $\tau$-algebra $\mathcal F$ of universe $F$ onto the full $\mathsf{FCA}$ $\mathbf B^{(\mathsf b)}$. 
\item $f$ is a t-homomorphism from the t-algebra $(A,\mathsf a,\varphi)_{\varphi\in F}$ onto the t-algebra $(B,\mathsf b,f^\bullet(\varphi))_{\varphi\in F}$, both of type $F$.

\end{enumerate}
\end{lemma}

\begin{proof} 
(1) By hypothesis $f^\mathbb N: \mathsf a\to \mathsf b$ is onto. 
Let $\varphi_0,\dots,\varphi_n\in F$ and let $\psi_0,\dots,\psi_n: \mathbf b\to B$ such that $\psi_i\circ f^\mathbb N=f\circ \varphi_i$. We have, for every $s\in\mathsf a$:
 \[
\begin{array}{llll}
&&f(q_n^\mathsf a(\varphi_0,\varphi_1,\dots,\varphi_n)(s))\\
&=&f(\varphi_0(s[\varphi_1(s),\dots,\varphi_n(s)]))\\
&=&\psi_0( f^\mathbb N(s[\varphi_1(s),\dots,\varphi_n(s)]))\\
&=&\psi_0( f^\mathbb N(s)[f(\varphi_1(s)),\dots,f(\varphi_n(s))])\\
&=&\psi_0( f^\mathbb N(s)[\varphi_1(f^\mathbb N(s)),\dots,\varphi_n(f^\mathbb N(s))])\\
&=&q_n^\mathsf b(\psi_0,\psi_1,\dots,\psi_n)(f^\mathbb N(s)).\\
\end{array}
\]
It follows that $F$  is a $\mathsf{FCA}$.  Moreover, $F$ is rich because every constant function trivially belongs to $F$.
%
%

(2) Given $\varphi\in F$, there exists a unique $\psi$ such that $\psi\circ f^\mathbb N=f\circ \varphi$. We denote $\psi$ by $f^\bullet(\varphi)$.
Then $f^\bullet$ can be proved to be a homomorphism by substituting $f^\bullet(\varphi_i)$ for $\psi_i$ in the chain of equalities of (i).
We now prove that $f^\bullet$ is onto. Let $\psi:\mathsf b\to B$ be any t-operation and let $ch: \mathcal P(A)\to A$ be a choice function (i.e., $ch(Y)\in Y$ for every nonempty subset $Y$ of $A$). We define $\varphi:\mathsf a \to A \in F$ as follows: $\varphi(s)=ch(f^{-1}(\psi(f^\mathbb N(s)))$.
Then $f^\bullet(\varphi)(f^\mathbb N s)=f(\varphi(s))= f(ch(f^{-1}(\psi(f^\mathbb N(s)))))=\psi(f^\mathbb N(s))$. Then $f^\bullet(\varphi)=\psi$.

(3)  By definition of $f^\bullet$, $f(\varphi(s))=f^\bullet(\varphi)( f^\mathbb N(s))$ for every $s\in \mathsf a$.
\end{proof}

\begin{theorem}\label{thm:giu}
 If  $H$ is a variety of clone $\tau$-algebras, then 
 
 (i) $H^\smalltriangledown$  is an Ft-variety of t-algebras. 
 
 (ii) $H^\blacktriangledown$  is an Et-variety of t-algebras. 
\end{theorem}
 
\begin{proof} (1) \emph{The classes $H^\smalltriangledown$ and $H^\blacktriangledown$ are closed under t-homomorphic image.}\\
Let $\mathbf A\in  H^\smalltriangledown$ and $f: \mathbf A\to \mathbf B$ be an onto t-homomorphism. We prove that $\mathbf B\in  H^\smalltriangledown$. 
Since $\mathbf A\in  H^\smalltriangledown$, then there exists a rich $\mathsf{FCA}_\tau$ $\mathcal G\in  H$ with value domain $\mathbf A$. Since $f$ is onto, then by Lemma \ref{lem:11} there exist a rich $\mathsf{FCA}_\tau$ $\mathcal F$ with value domain $\mathbf A$ and an onto homomorphism $f^\bullet$ from $\mathcal F$ onto $\mathbf B^{(\mathsf b)}$.  Since both $\mathcal F$ and $\mathcal G$ are rich and $\mathcal G\in H$, then $\mathcal F\cap \mathcal G\in H$ is rich. Then the image of $(f^\bullet)_{|\mathcal F\cap \mathcal G}$ is a rich $\mathsf{FCA}_\tau(\mathbf B)$ that belongs to $ H$. By definition of $H^\smalltriangledown$ this implies $\mathbf B\in  H^\smalltriangledown$. 
If $\mathbf A\in  H^\blacktriangledown$ the proof is similar.

(2) \emph{The classes $ H^\smalltriangledown$ and $ H^\blacktriangledown$ are closed under  t-subalgebra.}\\
If  $\mathbf A=(A,\mathsf a,\sigma^\mathbf A)\in  H^\smalltriangledown$, then there exists a rich $\mathsf{FCA}_\tau$ $\mathcal G\in  H$ with value domain $\mathbf A$. Let $\mathbf B=(B,\mathsf b,\sigma^\mathbf B)$ be a t-subalgebra of $\mathbf A$. 
We consider the subalgebra $\mathcal C$ of  $\mathcal G$ of all $\varphi:  \mathsf a \to A \in  G$ such that $\varphi_{|\mathsf b}$ is into $B$. $\mathcal C\in H$ and contains the semiconstant maps of values in $B$. The image of the  homomorphism from $\mathcal C$ into $\mathbf B^{(\mathsf b)}$, defined by $\varphi\in C\mapsto \varphi_{|\mathsf b}\in B^{\mathsf b}$, belongs to $H$ and it is a rich subalgebra of  $\mathbf B^{(\mathsf b)}$. Then $\mathbf B\in  H^\smalltriangledown$.
 If $\mathbf A\in  H^\blacktriangledown$ the proof is similar.

(3)  \emph{The classes $ H^\smalltriangledown$ and $ H^\blacktriangledown$ are closed under  t-product.}\\
Let $\mathbf B=\Pi_{i\in I}\mathbf A_i$ be the t-product of the t-algebras $\mathbf A_i\in  H^\smalltriangledown$ of trace $\mathsf a_i$, let $\mathcal F_i\in  H$ be a rich $\mathsf{FCA}_\tau(\mathbf A_i)$ ($i\in I$), and let $\mathsf b=\Pi_{i\in I}\mathsf a_i$. By Lemma \ref{sss}  $\Pi_{i\in I} \mathcal F_i$ can be embedded into the full  $\mathsf{FCA}$ $\mathbf B^{(\mathsf b)}$ and the image is rich if all $\mathcal F_i$ are rich. Then $\mathbf B\in   H^\smalltriangledown$.  If $\mathbf A_i\in  H^\blacktriangledown$ the proof is similar.

(4)  \emph{The class $ H^\blacktriangledown$ is closed under  t-expansion.}\\
Let $\mathbf A$ be a t-algebra of type $\tau$ and trace $\mathsf a$ such that the t-subalgebra $\mathbf A_{\bar s}$ generated by $s$ belongs to $ H^\blacktriangledown$, for every $s\in\mathsf a$. Then the term clone algebra $(\mathbf A_{\bar s})^\uparrow$ over $\mathbf A_{\bar s}$ belongs to $H$ for every $s\in \mathsf a$, because it is the minimal subalgebra of the full $\mathsf{FCA}_\tau$ $\mathbf A_{\bar s}^{(\mathsf a_{\bar s})}$ and by the hypothesis $\mathbf A_{\bar s}\in H^\blacktriangledown$ there exists a subalgebra of $\mathbf A_{\bar s}^{(\mathsf a_{\bar s})}$ belonging to $H$. By Proposition \ref{prop:c(a)} $\mathbf A^\uparrow$ is a subdirect product of the clone $\tau$-algebras $(\mathbf A_{\bar s})^\uparrow$. Then $\mathbf A^\uparrow\in  H$.
We conclude that $\mathbf A\in  H^\blacktriangledown$.

(5)  \emph{The class $ H^\smalltriangledown$ is closed under full t-expansion.}\\
Let $\mathbf A=\oplus^b_{i\in I} \mathbf A_{\upharpoonright\mathsf b_i}$ be a t-algebra of type $\tau$ and trace $\mathsf a$ such that the full t-subalgebra $\mathbf A_{\upharpoonright\mathsf b_i}$ belongs to $ H^\smalltriangledown$, for every $i\in I$. Then by definition of $ H^\smalltriangledown$ the polynomial clone algebra $(\mathbf A_{\upharpoonright\mathsf b_i})^{\!\Uparrow}$ belongs to $H$ for every $i\in I$. Since by Proposition \ref{prop:c(aa)} $\mathbf A^{\!\Uparrow}$ is a subdirect product of the clone $\tau$-algebras $(\mathbf A_{\upharpoonright\mathsf b_i})^{\!\Uparrow}$ through the map $t^\mathbf A\mapsto (t^\mathbf A_{|\mathsf b_i}:i\in I)$, then $\mathbf A^{\!\Uparrow}\in  H$.
We conclude that $\mathbf A\in  H^\smalltriangledown$.
\end{proof}

\begin{corollary}\label{cor:updown} Let $K$ be a class of t-algebras of type $\tau$. 
 If $K^{\smalltriangleup\smalltriangledown}=K$, then the following conditions are equivalent:
 \begin{enumerate}
\item $K$ is a t-variety;
\item $K$ is an Ft-variety;
\item $K^\smalltriangleup$ is a variety of $\mathsf{CA}_\tau$s.
\end{enumerate}
\end{corollary}

\begin{proof} (1 $\Rightarrow$ 3)  By Theorem \ref{thm:ggg}. (3 $\Rightarrow$ 2) By Theorem \ref{thm:giu}(i) and $K^{\smalltriangleup\smalltriangledown}=K$. (2 $\Rightarrow$ 1) Trivial.
\end{proof}

\begin{corollary}\label{cor:updown2}
If $K^{\smalltriangleup\blacktriangledown}=K$, then the following conditions are equivalent:
\begin{enumerate}
\item $K$ is a t-variety;
\item $K$ is an Ft-variety;
\item $K^\smalltriangleup$ is a variety of $\mathsf{CA}_\tau$s;
\item $K$ is an Et-variety.
\end{enumerate}
\end{corollary}

\begin{proof}
 By $K\subseteq K^{\smalltriangleup\smalltriangledown}\subseteq K^{\smalltriangleup\blacktriangledown}$ we derive $K^{\smalltriangleup\smalltriangledown}=K$. Then from Corollary \ref{cor:updown} it follows the equivalence of (1), (2) and (3).
(4) $\Rightarrow$ (3) By Theorem \ref{thm:ggg}. 
(3) $\Rightarrow$ (4) By hypothesis and Theorem \ref{thm:giu}(ii). 
\end{proof}

 \section{Birkhoff's Theorem for Et-varieties}\label{sec:gb}

Birkhoff's Theorem characterises the class $K$ of all models of an equational t-theory  of type $\tau$ in many different ways.
It generalises the classical Birkhoff's Theorem for varieties of algebras, but this is the subject of Section \ref{sec:newbir}.


\begin{lemma}\label{lem:ti} Let $\Sigma$ be a set of $\tau$-equations. Then $\mathrm{Mod}(\Sigma)=\{\mathbf A: \Sigma\subseteq \theta_\mathbf A\}$  is an Et-variety.
\end{lemma}

\begin{proof} 
Let $\mathbf A$ be a t-algebra of trace $\mathsf a$ such that the t-subalgebra $\mathbf A_{\bar s}$ generated by $s$ belongs to $\mathrm{Mod}(\Sigma)$ for every $s\in \mathsf a$. By Proposition \ref{prop:c(a)} $\theta_\mathbf A= \bigcap_{s\in \mathsf a} \theta_{\mathbf A_{\bar s}} \supseteq \Sigma$. Then $ \mathbf A\in \mathrm{Mod}(\Sigma)$.
The other closures are standard.
\end{proof}

\begin{theorem}\label{thm:bir1} {\rm (Birkhoff's Theorem for Et-varieties)} Let $K$ be a class of t-algebras of type $\tau$. Then the following conditions are equivalent:
\begin{enumerate}
\item $K$ is an Et-variety.
\item $K=\mathrm{Mod}(\theta_K)$.
\item $K=K^{\smalltriangleup\blacktriangledown}$ and $K^\smalltriangleup$ is a variety of clone $\tau$-algebras.
\item[(4)] $K=K^{\smalltriangleup\blacktriangledown}$ and $K$ is a t-variety.
\item[(5)] $K=K^{\smalltriangleup\blacktriangledown}$ and $K$ is an Ft-variety.
 \end{enumerate}
\end{theorem}

\begin{proof} The equivalence of (3), (4) and (5) follows from Corollary \ref{cor:updown2}.

(1 $\Rightarrow$ 2) Let $J$ be the set of all $\tau$-identities $e$ such that $e$ fails in some algebra $\mathbf A_e\in K$. 
The kernel of the unique homomorphism $f$  from $\mathcal T_\tau$ into  $\Pi_{e\in J} \mathbf A_e^\uparrow$
 is the equational t-theory $\theta_{K}$ of $K$.
Then $\mathcal T_K=\mathcal T_\tau/\theta_{K}$ is (up to isomorphism)  the minimal subalgebra of $\Pi_{e\in J} \mathbf A_e^\uparrow$.
By Lemma \ref{sss}  the clone $\tau$-algebra $\Pi_{e\in J} \mathbf A_e^\uparrow$ is isomorphic to a $\mathsf{FCA}_{\tau}$ $\mathcal F$ with value domain $\Pi_{e\in J}\mathbf A_e\in K$.  
 By  Proposition \ref{prop:zzz}(ii) the t-algebra $\mathcal F^\downarrow$ belongs to $\mathrm S_t\mathrm P_t(\Pi_{e\in J}\mathbf A_e)\subseteq K$. It follows that  the t-subalgebra $(\mathcal T_K)^\downarrow$ of  $\mathcal F^\downarrow$ belongs to $K$.


Let $\mathbf A\in \mathrm{Mod}(\theta_K)$ of trace $\mathsf a$. The unique homomorphism from $\mathcal T_\tau$ onto $\mathbf A^\uparrow$ can be factorised through  $\mathcal T_{K}$, because $\theta_K\subseteq \theta_\mathbf A$. Then $\mathbf A^\updownarrows\in K$, because it is a t-homomorphic image of $(\mathcal T_K)^\downarrow\in K$.
For every $s\in \mathsf a$, by Proposition \ref{prop:zzz}  the t-homomorphism $\bar s: \mathbf A^\updownarrows\to \mathbf A$, defined by $\bar s(t^\mathbf A)=t^\mathbf A(s)$, has the t-subalgebra $\mathbf A_{\bar s}$ generated by $s$ as image. Since  $\mathbf A_{\bar s}\in K$ for every $s\in\mathsf a$ and $K$ is closed under t-expansion, then  $\mathbf A\in K$. This concludes the proof that  $K=\mathrm{Mod}(\theta_K)$.

(2 $\Rightarrow$ 3) 
By Lemma \ref{lem:ti} $K$ is an Et-variety. Then by Theorem \ref{thm:ggg} $K^\smalltriangleup$ is a variety of clone $\tau$-algebras.
By definition we have $K\subseteq K^{\smalltriangleup\blacktriangledown}$. 
We now show that $K=K^{\smalltriangleup\blacktriangledown}$. Let $\mathbf A \in K^{\smalltriangleup\blacktriangledown}$, so that the term clone $\tau$-algebra $\mathbf A^\uparrow \in  K^\smalltriangleup$. Then there exist $\mathbf B\in K$ and an isomorphism $f:\mathbf A^\uparrow\to\mathbf B^\uparrow$ of clone $\tau$-algebras.   Since $\mathbf B\in K= \mathrm{Mod}(\theta_{K})$, then $\theta_{K}\subseteq \theta_{\mathbf B}\subseteq \theta_{\mathbf A}$.
 It follows that $\mathbf A\in K$.

%

(3 $\Rightarrow$ 1) By Corollary \ref{cor:updown2}.
\end{proof}

Question: Let $K$ be an Et-variety of t-algebras of type $\tau$ axiomatised by $\theta_K$. What is the axiomatisation of the variety $K^\smalltriangleup$ of clone $\tau$-algebras? The variety $K^\smalltriangleup$ is generated by the class of all full $\mathsf{FCA}_\tau$ with value domain $\mathbf A\in K$.
Then $\mathrm{Th}_{K^\smalltriangleup} = \mathrm{Tf}_K^X$, the f-hypertheory of $K$ (see Definition \ref{def:sat}).

Since $K^{\!\uparrow}\subseteq K^\smalltriangleup$, then  in general we have: $\mathrm{Th}_{K^\smalltriangleup}\subseteq \mathrm{Th}_{K^{\!\uparrow}}$.
For example, if $\tau=\emptyset$, then $K^\smalltriangleup=\mathsf{CA}_0$ (the class of all pure clone algebras) and $K^{\!\uparrow}=\mathbb I \{\mathcal P\}$, where $\mathcal P=(\mathbb N, q_n^\mathcal P,\e_i^\mathcal P)$ is the minimal pure  clone algebra defined in Example \ref{exa:proj}.
The identity $y(y(x_{11},x_{12},\e_3,\dots),y(x_{21},x_{22},\e_3,\dots),\e_3,\dots) = y(x_{11},x_{22},\e_3,\dots)$ belongs to $\mathrm{Th}_{K^{\!\uparrow}}\setminus \mathrm{Th}_{K^\smalltriangleup}$.

\section{Birkhoff's Theorem for Ft-varieties}\label{sec:gb2}
In this section we characterise the class $K$ of all models of a w-hypertheory  of type $\tau$ in different ways: 
\begin{enumerate}
\item $K$ is an Ft-variety of t-algebras of type $\tau$.
\item $K^\smalltriangleup$ is a variety of clone $\tau$-algebras and $K=K^{\smalltriangleup\smalltriangledown}$.
\item $K$ is a t-variety and $K=K^{\smalltriangleup\smalltriangledown}$.
\end{enumerate}

\noindent If $\Gamma$ is a set of hyperidentities,  $\mathrm{Mod}_w\Gamma=\{\mathbf A: \forall t_1=t_2\in\Gamma.\ \mathbf A\models_w t_1=t_2\}$. 
If $K$ is a class of t-algebras, then we denote by $\mathcal T_{K,X}$ the quotient of $\mathcal T_\tau(X)$ by the equational w-hypertheory $\mathrm{Tw}^X_K$.

\begin{lemma}\label{mmm} Let $K$ be a class of t-algebras of the same type and $\Gamma$ be a set of $\tau$-hyperidentities. Then, we have:
\begin{enumerate}
\item $\mathrm{Mod}_w(\mathrm{Tw}^X_K)= \mathrm{Mod}_w(\mathrm{Tw}^Y_K)$ for all $X,Y$ of infinite cardinality.
\item  $\mathrm{Mod}_w(\Gamma)$  is an Ft-variety.
\end{enumerate}
\end{lemma}

\begin{theorem}\label{thm:bir2} {\rm (Birkhoff's Theorem for Ft-varieties)} Let $K$ be a class of t-algebras of type $\tau$. Then the following conditions are equivalent:
\begin{itemize}
\item[(i)]  $K$ is an Ft-variety.
\item[(ii)] $K=\mathrm{Mod}_w(\mathrm{Tw}^X_K)$, for every infinite set $X$.
\item[(iii)] $K^\smalltriangleup$ is a variety of clone $\tau$-algebras and $K=K^{\smalltriangleup\smalltriangledown}$.
\item[(iv)] $K$ is a t-variety and $K=K^{\smalltriangleup\smalltriangledown}$.
\end{itemize}
\end{theorem}

\begin{proof} (i)  $\Rightarrow$ (ii) 
 Let $X$ be an infinite set and $I$ be the set of all hyperidentities $e\in T_\tau(X)^2$ such that $e$ weakly fails in some algebra $\mathbf A_e\in K$. In other words, if $e\equiv t_1=t_2$, then $\alpha_e(t_1)\neq \alpha_e(t_2)$  for some weak homomorphism $\alpha_e: \mathcal T_\tau(X) \to \mathbf A_e^{\!\Uparrow}$. 

Consider the unique homomorphism $f$  from $\mathcal T_\tau(X)$ into  $\Pi_{e\in I}  \mathbf A_e^{\!\Uparrow}$ such that $f(t) =(\alpha_e(t):e\in I)$. 
The kernel of $f$ is the equational w-hypertheory $\mathrm{Tw}^X_K$ of $K$.
Then $\mathcal T_{K,X}$ is up to isomorphism  a subalgebra of $\Pi_{e\in I}  \mathbf A_e^{\!\Uparrow}$.
By Lemma \ref{sss}  the clone $\tau$-algebra $\Pi_{e\in I}  \mathbf A_e^{\!\Uparrow}$ is isomorphic to a $\mathsf{FCA}_{\tau}$ $\mathcal F$ with value domain $\Pi_{e\in I}\mathbf A_e\in K$.  
 By  Proposition \ref{prop:zzz}(ii) the t-algebra $\mathcal F^\downarrow$ belongs to $\mathbb S\mathbb P(\Pi_{e\in I}\mathbf A_e)\subseteq K$. It follows that  the t-subalgebra $(\mathcal T_{K,X})^\downarrow$ of the t-algebra $\mathcal F^\downarrow$ belongs to $K$.

Let $\mathbf A\in \mathrm{Mod}_w(\mathrm{Tw}^X_K)$ of trace $\mathsf a$ and $Y$ be an infinite set, whose cardinality is greater than the cardinality of  $A$. Let $\alpha:\mathcal T_\tau(Y) \to \mathbf A^{\!\Uparrow}$ be an onto weak homomorphism.
By Lemma \ref{mmm} $\mathbf A\in \mathrm{Mod}_w(\mathrm{Tw}^Y_K)$.
 Since $\mathbf A^{\!\Uparrow}$ is a homomorphic image of  $\mathcal T_{K,Y}$, then $\mathbf A^{\!\Uparrow\downarrow}$ is a homomorphic image of  $(\mathcal  T_{K,Y})^\downarrow\in K$.
Then $\mathbf A^{\!\Uparrow\downarrow}\in K$.
The image of the t-homomorphism $\bar r: \mathbf A^{\!\Uparrow\downarrow}\to \mathbf A$, defined by $\bar r(\varphi)=\varphi(r)$ for every $\varphi\in \mathbf A^{\!\Uparrow\downarrow}$, is the full subalgebra of $\mathbf A$ of basic trace $[r]_\mathbb N^A$. Then the closure under $F_t$ ensures that
 $\mathbf A\in K$. 

(ii)  $\Rightarrow$ (iii) By Lemma \ref{mmm} $K$ is an Ft-variety. Then by Theorem \ref{thm:ggg} $K^\smalltriangleup$ is a variety of clone $\tau$-algebras.
By definition we have $K\subseteq K^{\smalltriangleup\smalltriangledown}$. 
We now show that $K=K^{\smalltriangleup\smalltriangledown}$. Let $\mathbf A \in K^{\smalltriangleup\smalltriangledown}$ of trace $\mathsf a$, so that the polynomial clone $\tau$-algebra $\mathbf A^{\!\Uparrow} \in  K^\smalltriangleup$. Then there exist $\mathbf B\in K$ of trace $\mathsf b$, and an embedding $f: \mathbf A^{\!\Uparrow}\to\mathbf B^{\!(\mathsf b)}$ of clone $\tau$-algebras. Then, $f(\varphi)$ is semiconstant, for every constant t-operation $\varphi:\mathsf a\to A$.
We now show that $\mathrm{Tw}^X_\mathbf B\subseteq \mathrm{Tw}^X_\mathbf A$. Let $\mathbf B\models_w u_1 =u_2$ for a hyperidentity $u_1=u_2$. For every weak homomorphism $\alpha: \mathcal T_\tau(X)\to \mathbf A^{\!\Uparrow}$, $f\circ \alpha: \mathcal T_\tau(X)\to \mathbf B^{\!(\mathsf b)}$ is  a strong homomorphism. By Proposition \ref{prop:ws} $\mathbf B\models_w u_1 =u_2$ iff $\mathbf B\models_s u_1 =u_2$.
Since $f$ is injective, then $\alpha(u_1)=\alpha(u_2)$ for every $\alpha$, i.e., $\mathbf A\models_w u_1 =u_2$.
  Since $\mathbf B\in K= \mathrm{Mod}_w(\mathrm{Tw}_K^X)$, then $\mathrm{Tw}_K^X\subseteq \mathrm{Tw}^X_\mathbf B\subseteq \mathrm{Tw}^X_\mathbf A$.
 It follows that $\mathbf A\in K$.

(iii)  $\Rightarrow$ (i)  and (iii) $\Leftrightarrow$ (iv) follow from Corollary \ref{cor:updown}.
\end{proof}

\section{A new Birkhoff's theorem for varieties of algebras}\label{sec:newbir}
In this section we show a new version of Birkhoff HSP theorem for varieties of algebras.
We prove that a class $H$ of algebras is a variety iff the class $H^\star$ of t-algebras, obtained by gluing together the algebras in $H$,  satisfies:  $(H^\star)^\smalltriangleup$ is a variety of clone algebras and $H^\star=(H^\star)^{\smalltriangleup\blacktriangledown}$.

\subsection{t-Algebras from  algebras} \label{sec:tafromalg}
In this section we study the t-algebras obtained by gluing together a family of algebras.

Recall that, if $\rho=(\rho_n: n\geq 0)$ is a finitary type, then $\rho^\star=\bigcup_{n\geq 0}\rho_n$.


Let $\rho$ be a finitary type and $\mathbf S=(S,\sigma^\mathbf S)_{\sigma\in\rho}$ be a $\rho$-algebra. For every trace $\mathsf a$ on $S$,
we denote by $\mathbf S^\top_\mathsf a$ the t-algebra of type $\rho^\star$ and trace $\mathsf a$ such that $\sigma^{\mathbf S^\top_\mathsf a}= (\sigma^\mathbf S)^\top_\mathsf a$ (see Example \ref{exa:extfinalg}).

\begin{definition} 
\begin{itemize}
\item[(i)]  A t-algebra $\mathbf A$ of type $\tau$ is \emph{finite-dimensional} if $\mathbf A^\uparrow$ is a finite-dimensional clone $\tau$-algebra.
\item[(ii)]  A t-algebra $\mathbf A$ of type $\tau$ is \emph{$\rho$-dimensional}, where $\rho$ is a finitary type, if $\tau=\rho^\star$ and $\mathrm{dim}(\sigma^\mathbf A)\leq n$ for every $\sigma\in \rho_n$.
\end{itemize}
\end{definition}


\begin{lemma}\label{lem:otto} Let $\mathbf A$ be a t-algebra of type $\tau$ and trace $\mathsf a$, and $\{\mathsf b_i\}_{i\in I}$ be the basic t-partition of $\mathsf a$. Then $\mathbf A$ is finite-dimensional if and only if
 $\mathbf A = \oplus^b_{i\in I} (\mathbf A_i)_{\mathsf b_i}^\top$ and the algebras $\mathbf A_i$ ($i\in I$) are $\rho$-algebras for some finitary type  $\rho$.
\end{lemma}

\begin{proof} ($\Rightarrow$) We define the finitary type $\rho$ as follows: for every $\sigma\in\tau$,  $\sigma\in\rho_n$ iff $\mathrm{dim}(\sigma^\mathbf A)=n$. Moreover, for every $i\in I$, we define the $\rho$-algebra
$\mathbf A_i = (A, \sigma^{\mathbf A_i})_{\sigma\in \rho}$ as follows for every $\sigma\in\rho_n$: $\sigma^{\mathbf A_i}(a_1,\dots,a_n)=\sigma^\mathbf A(s[a_1,\dots,a_n])$ for every $a_1,\dots,a_n\in A$ and $s\in\mathsf b_i$.

($\Leftarrow$) $(\mathbf A_i)_{\mathsf b_i}^\top$ is finite-dimensional for every $i\in I$.
\end{proof}

\begin{corollary} In the hypothesis of Lemma \ref{lem:otto}, 
 $\mathbf A$ is $\rho$-dimensional if and only if
 $\mathbf A = \oplus^b_{i\in I} (\mathbf A_i)_{\mathsf b_i}^\top$ for some $\rho$-algebras $\mathbf A_i$ ($i\in I$).
\end{corollary}

Let  $F_\rho(I)$ be  the set of $\rho$-terms in the variables $I=\{v_1,v_2,\dots\}$. We define a map $(-)^\star: F_\rho(I)\to T_{\rho^\star}$ as follows:
\begin{enumerate}
\item $(v_i)^\star= \e_i$ for every $v_i\in V$;
\item $\sigma(t_1,\dots,t_n)^\star= \sigma(t_1^\star,\dots,t_n^\star,\e_{n+1},\e_{n+2},\dots)$, for every $\sigma\in\rho_n$.
\end{enumerate}
The map $($-$)^\star$ is not onto. For example, if $\sigma \in \rho_2$, then $\sigma(\e_3,\e_2,\e_1,\e_4,\e_5,\dots)$ is not in the image of $(-)^\star$.

We define the inverse translation $(-)^\bullet: T_{\rho^\star}\to F_\rho(I)$ as follows:

\begin{enumerate}
\item $(\e_i)^\bullet= v_i$;
\item For every $\sigma\in\rho_n$,
$$\sigma(t_1,\dots,t_k,\e_{k+1},\e_{k+2},\dots)^\bullet=\begin{cases}\sigma(t_1^\bullet,\dots,t_k^\bullet,v_{k+1},\dots,v_n)&\text{if $k\leq n$}\\
\sigma(t_1,\dots,t_n) &\text{if $k> n$}
\end{cases}$$
\end{enumerate}

\begin{lemma}\label{lem:upterms1} Let $\mathbf A$ be a $\rho$-dimensional t-algebra. Then we have, for every $\rho^\star$-term $t$ and $\rho$-term $p$:
\begin{itemize}
\item[(i)] $\mathbf A\models (t^\bullet)^\star  = t$.

\item[(ii)] $(p^\star)^\bullet\equiv p$, where $\equiv$ is the syntactical identity.
\end{itemize}
\end{lemma}

\begin{proof}
(i)-(ii)  The proof is by induction over the complexity of $t$ and $p$, respectively.  
\end{proof}

\begin{lemma}\label{lem:upterms2} Let $\mathbf A = \oplus^b_{i\in I} (\mathbf A_i)_{\mathsf b_i}^\top$ be a $\rho$-dimensional t-algebra of trace $\mathsf a$, $t,u\in T_{\rho^\star}(\e_1,\dots,\e_n)$ be $\rho^\star$-terms and $p,z\in F_\rho(v_1,\dots,v_n)$ be $\rho$-terms.  Then we have: 

\begin{itemize}
\item[(i)] $(p^\star)^\mathbf A(s)=p^{\mathbf A_i}(s_1,\dots,s_n)$ and $t^\mathbf A(s)= (t^\bullet)^{\mathbf A_i}(s_1,\dots,s_n)$, for every  $s\in\mathsf b_i$.

\item[(ii)]  $\mathbf A\models p^\star=z^\star$ iff $\mathbf A_i \models p=z$ for every $i\in I$.

\item[(iii)] $\mathbf A\models t=u $ iff $\mathbf A\models (t^\bullet)^\star=(u^\bullet)^\star$ iff  $\mathbf A_i \models t^\bullet=u^\bullet$ for every $i\in I$.

\item[(iv)] $\theta_\mathbf A = \{t=u :  t^\bullet =u^\bullet \in \bigcap_{i\in I}\mathrm{Th}_{\mathbf A_i}\}$.
\end{itemize}
\end{lemma}

\begin{proof}
(i)  We have $\sigma^{\mathbf A_i}(a_1,\dots,a_n)=\sigma^\mathbf A(s[a_1,\dots,a_n])$ for every thread $s\in \mathsf b_i$ and $\sigma\in \rho_n$.
Then the proof is by induction over the complexity of $p$ and $t$, respectively. 

(ii) By (i).

(iii) By (ii) and Lemma \ref{lem:upterms1}(i).

(iv) follows from (iii). 
\end{proof}

\begin{example}
Let $\rho$ be a finitary type, $\sigma\in \rho_2$ and $t=\sigma(\e_1,\e_2,\sigma)$ be a $\rho^\star$-term.  Then $t^\bullet = \sigma(\e_1,\e_2)$ and $(t^\bullet)^\star = \sigma$. Every $\rho$-dimensional t-algebra satisfies the identity $\sigma(\e_1,\e_2,\sigma)=\sigma$.
\end{example}

\begin{definition}
 A $\rho^\star$-identity $(t^\bullet)^\star  = t$ ($t\in T_{\rho^\star}$) is called \emph{$\rho$-structural} if  $t$ is not in the image of the map $(-)^\star$.
\end{definition}

The set of all structural $\rho^\star$-identities associated with the finitary type $\rho$ will be denoted by $\mathrm{Str}(\rho)$.

\begin{example}\label{exa:stru}
 The identities $\sigma=\sigma(\e_1,\dots,\e_{k-1},\e_{k+1})$ ($k>n$, $\sigma\in\rho_n$), specifying that $\mathrm{dim}(\sigma^\mathbf A)\leq n$, are $\rho$-structural. 
\end{example}

\begin{proposition}\label{prop:str}
 The class of all $\rho$-dimensional t-algebras is an Et-variety of type $\rho^\star$ axiomatised by $\mathrm{Str}(\rho)$.
\end{proposition}

\begin{proof}
 By Lemma \ref{lem:upterms1}(i) every $\rho$-dimensional t-algebra satisfies $\mathrm{Str}(\rho)$. The converse follows from Example \ref{exa:stru}.
\end{proof}

\subsection{Birkhoff for algebras from Birkhoff for t-algebras}\label{sec:bfafbfta}
In this section we provide the proof of the new version of Birkhoff's theorem.

As a matter of notation, if $\mathbf A= \oplus^b_{i\in I} (\mathbf A_i)_{\mathsf b_i}^\top$ is $\rho$-dimensional, then we put $\mathrm{Th}_{\mathbf A,\rho}=\bigcap_{i\in I}\mathrm{Th}_{\mathbf A_i}$. $\mathrm{Th}_{\mathbf A,\rho}$ is a set of $\rho$-identities.

\begin{definition}\label{def:t-extH}
  Let $H$ be a class of algebras of finitary type $\rho$. The \emph{t-algebra extension of $H$} is the class
 $H^\star$  of all $\rho$-dimensional t-algebras $\mathbf A= \oplus^b_{i\in I} (\mathbf A_i)_{\mathsf b_i}^\top$ such that $\mathbf A_i\in H$ for every  $i\in I$.
\end{definition}

We have that $\mathrm{Th}_H=\bigcap_{\mathbf A\in H^\star}\mathrm{Th}_{\mathbf A,\rho}$. 

 If $E$ is a set of $\rho$-identities, then $E^\star=\{t_1^\star=t_2^\star : t_1=t_2\in E\}$ is the corresponding set of $\rho^\star$-identities.

\begin{lemma}\label{lem:H*} 
\begin{itemize}
\item[(i)] $\theta_{H^\star}=\{t=u: t^\bullet=u^\bullet \in \mathrm{Th}_H\}$.
\item[(ii)]  $(\mathrm{Th}_H)^\star \cup \mathrm{Str}(\rho) \subseteq \theta_{H^\star}$.
\item[(iii)] $\mathrm{Mod}((\mathrm{Th}_H)^\star \cup \mathrm{Str}(\rho)) =\mathrm{Mod}(\theta_{H^\star})$.
\end{itemize}

\end{lemma}

\begin{proof} (i) By Lemma \ref{lem:upterms2}(iii).

(ii) By Proposition \ref{prop:str} $\mathrm{Str}(\rho) \subseteq \theta_{H^\star}$. Let now $p=z\in \mathrm{Th}_H$.  By (i) $p^\star = z^\star \in \theta_{H^\star}$ iff $(p^\star)^\bullet=(z^\star)^\bullet \in \mathrm{Th}_H$. The conclusion follows from Lemma \ref{lem:upterms1}(ii) and the hypothesis: $(p^\star)^\bullet =p= z=(z^\star)^\bullet \in \mathrm{Th}_H$.

(iii) Let $\mathbf A\models (\mathrm{Th}_H)^\star \cup \mathrm{Str}(\rho)$.  By Proposition \ref{prop:str} $\mathbf A$ is $\rho$-dimensional.   Let $t=u\in \theta_{H^\star}$. Then  by (i) $t^\bullet=u^\bullet\in \mathrm{Th}_H$, so that  by the hypothesis $\mathbf A\models (\mathrm{Th}_H)^\star$ we derive  $\mathbf A\models (t^\bullet)^\star=(u^\bullet)^\star$. Since $\mathbf A$ is $\rho$-dimensional, we conclude by applying Lemma \ref{lem:upterms1}(i):  $\mathbf A\models t=(t^\bullet)^\star=(u^\bullet)^\star=u$.
%
\end{proof}

\begin{lemma}\label{lem:*} Let $H$ be a class of $\rho$-algebras. The following conditions are equivalent:
\begin{enumerate}
\item $H = \mathrm{Mod}(\mathrm{Th}_H)$;
\item $H^\star =\mathrm{Mod}(\theta_{H^\star})$.
\end{enumerate}
\end{lemma}

\begin{proof} (1) $\Rightarrow$ (2) 
Let $\mathbf A\in \mathrm{Mod}(\theta_{H^\star})$ be a t-algebra. By Lemma \ref{lem:H*}(ii) $\mathbf A\models \mathrm{Str}(\rho)$. Then
$\mathbf A$ is $\rho$-dimensional and $\mathbf A= \oplus^b_{i\in I} (\mathbf A_i)_{\mathsf b_i}^\top$. By  applying again Lemma \ref{lem:H*}(ii) we get $\mathbf A\models (\mathrm{Th}_H)^\star$. Then by Lemma \ref{lem:upterms2}(ii) $\mathbf A_i\models \mathrm{Th}_H$ for every $i\in I$, so that $\mathbf A\in H^\star$.
 
 (2) $\Rightarrow$ (1)  Let $\mathbf S\in \mathrm{Mod}(\mathrm{Th}_H)$ be a $\rho$-algebra and $\mathsf d$ be a basic trace on $S$. Then the t-algebra $\mathbf S_{\mathsf d}^\top$ over $\mathbf S$  is $\rho$-dimensional.  By  Lemma \ref{lem:upterms2}(ii) and Proposition \ref{prop:str} $\mathbf S_{\mathsf d}^\top$ satisfies  $(\mathrm{Th}_H)^\star \cup \mathrm{Str}(\rho)$. Thus $\mathbf S_\mathsf d^\top\in H^\star$ and $\mathbf S\in H$.
\end{proof}

%
%
%
%
%

Let $\bar x=x_1,\dots,x_k$, $\bar\e= \e_{m+1},\dots,\e_{m+k}$  and $t=t(\bar x)$ be  a $\rho^\star$-hyperterm. We define $t[\bar\e/\bar x]$ as follows:
 \begin{itemize}
\item $\e_j[\bar\e/\bar x]= \e_j$;
\item $\sigma(t_1,\dots,t_j,\e_{j+1},\dots)[\bar\e/\bar x]= \sigma(t_1[\bar\e/\bar x],\dots,t_j[\bar\e/\bar x],\e_{j+1},\dots)$;
\item $x_i(t_1,\dots,t_j,\e_{j+1},\dots)[\bar\e/\bar x]= \e_{m+i}$.
\end{itemize}

\begin{lemma}\label{lem:ninuzzo} Let  $\rho$ be a finitary type, $\bar x=x_1,\dots,x_k$, $t=t(\bar x)$ be  a $\rho^\star$-hyperterm and  $K$ be a  class of $\rho$-dimensional  t-algebras. Then there exists a natural number $l(t)$ such that, for every $m\geq l(t)$, $\mathbf A\in K$ of trace $\mathsf a$,  $s\in \mathsf a$ and $\alpha: \bar x\to A$, we have 
$$\alpha^\star(t)(s)=t[\bar\e/\bar x]^\mathbf A(s[s_1,\dots,s_m,\alpha(x_1),\dots,\alpha(x_k)]),$$
where $\alpha^\star:\mathcal T_\tau(\bar x)\to \mathbf A^{\!\Uparrow}$ is the weak homomorphism extending $\alpha$. 
\end{lemma}

\begin{proof} Let $\mathbf A = \oplus^b_{i\in I} (\mathbf A_i)_{\mathsf b_i}^\top\in K$ be a $\rho$-dimensional t-algebra of trace $\mathsf a$, where $\mathbf A_i$ is a $\rho$-algebra. Let $1\leq i\leq k$.

$t\equiv \e_i$: $m=l(\e_i)=i$.

$t\equiv x_i(t_1,\dots,t_j,\e_{j+1},\dots)$: We put $m=l(t)=0$. Then $\alpha^\star(t)(s)=\alpha(x_i)$ for every $s\in\mathsf a$, because $\alpha^\star(x_i)$ is the constant function of value $\alpha(x_i)$. Moreover, 
$t(\bar\e/\bar x)= \e_i$ and  $\e_i^\mathbf A(\alpha(x_1),\dots,\alpha(x_k),s_{k+1},\dots)=\alpha(x_i)$.

$t\equiv \sigma(t_1,\dots,t_j,t_{j+1},\dots)$ with $\sigma\in\rho_n$, $t_j\neq \e_j$ and $t_k=\e_k$ for $k>j$:\\
 By applying the induction hypothesis, let $l(t_i)$ ($i=1,\dots,n$) be the natural number specified in this lemma.
 We define $l(t)=1+\mathrm{max}\{n,l(t_1),\dots,l(t_n)\}$. Let $m\geq l(t)$, $\bar \e=\e_{m+1},\dots,\e_{m+k}$ and $f_i:A^m\to A$ ($i\in I$) such that $(\sigma^\mathbf A)_{|\mathsf b_i}=(f_i)^\top_{\mathsf b_i}$. Let $s\in\mathsf b_i$ and $r= s[s_1,\dots,s_m,\alpha(x_1),\dots,\alpha(x_k)]\in\mathsf b_i$.
 Then we have:
 \[
\begin{array}{lll}
 \alpha^\star(t)(s) & =  &  \alpha^\star(\sigma(t_1,\dots,t_j,\e_{j+1},\dots))(s) \\
 & =  &  \alpha^\star(q_m^\mathcal T(\sigma, t_1,\dots,t_j,\e_{j+1},\dots,\e_m))(s) \\
  &  = & q_m^\mathsf a(\sigma^\mathbf A, \alpha^\star(t_1),\dots,\alpha^\star(t_j),\alpha^\star(\e_{j+1}),\dots,\alpha^\star(\e_m))(s)  \\
    &  = & \sigma^\mathbf A(s[\alpha^\star(t_1)(s),\dots,\alpha^\star(t_j)(s), s_{j+1},\dots,s_m])  \\

  & =  &   f_i(\alpha^\star(t_1)(s),\dots,\alpha^\star(t_j)(s), s_{j+1},\dots,s_m)\\
    & =  &f_i( t_1[\bar\e/\bar x]^\mathbf A(r),\dots, t_j[\bar\e/\bar x]^\mathbf A(r),\e_j[\bar\e/\bar x]^\mathbf A(r),\dots,\e_m[\bar\e/\bar x]^\mathbf A(r)) \\
& =  &\sigma^\mathbf A(t_1[\bar\e/\bar x]^\mathbf A(r),\dots,t_j[\bar\e/\bar x]^\mathbf A(r),\dots) \\
& =  & t[\bar\e/\bar x]^\mathbf A(r).\qedhere
\end{array}
\]
\end{proof}

\begin{corollary}\label{cor:wt}
 In the hypotheses of Lemma \ref{lem:ninuzzo}, let $\bar\e=\e_{m+1},\dots,\e_{m+k}$ for  $m>l(t_1),l(t_2)$. Then $\mathbf A\models_w t_1(\bar x)=t_2(\bar x)$ iff $\mathbf A\models t_1[\bar\e/\bar x]=t_2[\bar\e/\bar x]$.

\end{corollary}

\begin{corollary}\label{cor:tetvar}
Let  $\rho$ be a finitary type and  $K$ be a   class of $\rho$-dimensional t-algebras. Then  $K$ is an Ft-variety iff 
it is an Et-variety.
\end{corollary}

\begin{proof}
By Corollary \ref{cor:wt} and Theorems \ref{thm:bir1}, \ref{thm:bir2}. 
\end{proof}

\begin{theorem} {\rm (Birkhoff's Theorem for Algebras)}\label{thm:birkhoffforalgebras} Let $\rho$ be a finitary type. The following conditions are equivalent for a class $H$ of  $\rho$-algebras:
\begin{enumerate}
\item $H$ is a variety; 
\item $H=\mathrm{Mod}(\mathrm{Th}_H)$;
\item  $H^\star$ is an Et-variety; 
\item $H^\star=\mathrm{Mod}(\theta_{H^\star})$;
\item $(H^\star)^\smalltriangleup$ is a variety of clone $\rho^\star$-algebras and $H^\star=(H^\star)^{\smalltriangleup\blacktriangledown}$;
\item  $H^\star$ is an Ft-variety;
\item $H^\star=\mathrm{Mod}_w(\mathrm{Tw}^X_{H^\star})$, for an infinite set $X$;
\item $(H^\star)^\smalltriangleup$ is a variety of clone $\rho^\star$-algebras and $H^\star=(H^\star)^{\smalltriangleup\smalltriangledown}$;
\item  $H^\star$ is a t-variety.
\end{enumerate}
\end{theorem}

\begin{proof} (1) $\Leftrightarrow$ (2) is the classical $\mathrm{HSP}$ Birkhoff Theorem (see \cite[Theorem 11.9]{BS81}).
The equivalence of (3), (4) and (5) follows from Theorem \ref{thm:bir1}.
(2) $\Leftrightarrow$ (4)  is shown in Lemma \ref{lem:*}.
 By Theorem \ref{thm:bir2} items (6),(7) and (8) are equivalent. By Corollary \ref{cor:tetvar} we get the equivalence of (6) and (3).\\
 (9) $\Rightarrow$ (1) We show that $H$ is closed under $\mathbb H$, $\mathbb S$ and $\mathbb P$.
 \begin{itemize}
\item  Let $\mathbf S\in H$ be a $\rho$-algebra and $\mathbf T$ be a subalgebra of $\mathbf S$. Since $\mathbf S^\top\in H^\star$ and $\mathbf T^\top$ is the t-subalgebra of $\mathbf S^\top$ determined by the trace $T^\mathbb N \subseteq S^\mathbb N$ on $T$, then $\mathbf T^\top\in H^\star$. By definition of $H^\star$ we obtain that $\mathbf T\in H$.

\item  Let $\mathbf S\in H$ and $f: \mathbf S\to \mathbf W$ be an onto homomorphism of algebras. Since $f$ is also an onto t-homomorphism from $\mathbf S^\top$ into $\mathbf W^\top$, then $\mathbf W^\top\in H^\star$ and we conclude $\mathbf W\in H$.

\item Let $\mathbf S_i\in H$ ($i\in I$) be a family of algebras and $\mathbf A=\prod_{i\in I}\mathbf S_i^\top$.  Since $\mathbf S_i^\top\in H^\star$, then by the hypothesis we obtain $\mathbf A\in H^\star$. The universe of $\mathbf A$ is the set $\prod_{i\in I} S_i$ and its trace is $\mathsf a = \prod_{i\in I} S_i^\mathbb N = (\prod_{i\in I} S_i)^\mathbb N$. It remains to show that $\sigma^\mathbf A = \prod_{i\in I} \sigma^{\mathbf S_i^\top}$.
 This is the case, because $\sigma^\mathbf A(s)= (\sigma^{\mathbf S_i^\top}(s^i) :i\in I)$ for every $s= (s^i : i\in I)\in\mathsf a$.
 Then $\prod_{i\in I}\mathbf S_i^\top = (\prod_{i\in I}\mathbf S_i)^\top$ and by definition of $H^\star$ we conclude $\prod_{i\in I}\mathbf S_i\in H$.
 \qedhere
\end{itemize}
\end{proof}

\section{Free t-algebras}\label{sec:free}

In this section we characterise the t-algebras satisfying $\mathbf A\cong\mathbf A^\updownarrows$.

In the following $K$ is a class of t-algebras of type $\tau$ and $\mathbf A$ is a t-algebra of the same type and trace $\mathsf a$.

\subsection{The term clone algebra as initial algebra}
We say that \emph{a clone $\tau$-algebra $\mathcal C$ has the initial property for a class $X$ of clone $\tau$-algebras} if $\mathcal C$ is minimal (see Definition \ref{def:min}) and, for every $\mathcal D\in  X$,
there exists a unique homomorphism from $\mathcal C$ into $\mathcal D$.

We denote by $\mathrm V_E(\mathbf A)$ the Et-variety generated by the t-algebra $\mathbf A$, i.e., $\mathrm V_E(\mathbf A) = \mathrm{Mod}(\theta_\mathbf A)$.

\begin{proposition}\label{prop:tre} The following conditions are equivalent:
\begin{enumerate}
\item $\mathbf A^\uparrow$ has the initial property for $K^\smalltriangleup$; 
\item $\theta_\mathbf A\subseteq \theta_K$;
\item $K\subseteq \mathrm V_E(\mathbf A)$.
\end{enumerate}
If $\mathbf A\in K$, then the above three conditions are equivalent to 
\begin{enumerate}
\item[(4)] $\mathrm V_E(K)= \mathrm V_E(\mathbf A)$ (i.e., $\theta_\mathbf A=\theta_K$).
\end{enumerate}
\end{proposition}

\begin{proof}
  $\mathbf B\in \mathrm V_E(\mathbf A)$ iff  (by Theorem \ref{thm:bir1}) $\theta_\mathbf A\subseteq \theta_\mathbf B$ iff  $t^\mathbf A\mapsto t^\mathbf B$ is a homomorphism from $\mathbf A^\uparrow$ onto $\mathbf B^\uparrow$.
\end{proof}

\subsection{Free t-algebras and minimal clone algebras}

\begin{definition}\label{def:free} 
  We say that \emph{$\mathbf A$ has the universal mapping property for $K$} if there exists $s\in \mathsf a$ satisfying the following properties:
\begin{enumerate}
\item $\mathbf A=\mathbf A_{\bar s}$ is generated by $s$. 
\item For every $\mathbf B\in K$ of trace $\mathsf b$ and $r\in  \mathsf b$, there is a unique t-homomorphism $\alpha:\mathbf A\to\mathbf B$ such that $\alpha^\mathbb N(s)=r$.
\end{enumerate}
\end{definition}

Note that, if the t-algebra $\mathbf A$ of trace $\mathsf a$ has the universal mapping property for $K$, then the trace $\mathsf a$ is basic.

\begin{proposition}\label{prop:quattro} The following conditions are equivalent for a t-algebra $\mathbf A\in K$ of trace $\mathsf a$:
\begin{enumerate}
\item $\mathbf A$ has the universal mapping property for $K$;
\item $\mathbf A\cong\mathbf A^\updownarrows$ and $\mathbf A^\uparrow$ has  the initial property for $K^\smalltriangleup$.
\end{enumerate}
 $\mathbf A$ will be called the free $K$-t-algebra. 
\end{proposition}

\begin{proof} (1) $\Rightarrow$ (2) By hypothesis there is $s\in \mathsf a$ satisfying the hypotheses of Definition \ref{def:free}. 
The map $[t]_{\theta_\mathbf A}\mapsto t^\mathbf A(s)$ is a  t-isomorphism from $(\mathcal T_\tau/\theta_\mathbf A)^\downarrow$ onto $\mathbf A$. 
Since $\mathbf A^\uparrow\cong \mathcal T_\tau/\theta_\mathbf A$, then 
 it follows that $\mathbf A^\updownarrows\cong (\mathcal T_\tau/\theta_\mathbf A)^\downarrow \cong \mathbf A$.


We now prove that $\mathbf A^\uparrow$ has  the initial property for $K^\smalltriangleup$. It is sufficient to show that $\theta_\mathbf A\subseteq \theta_\mathbf B$ for every $\mathbf B\in K$ of trace $\mathsf b$. For every $r\in \mathsf b$, there exists a unique onto t-homomorphism  $f_r:   \mathbf A\to \mathbf B_{\bar r}$  such that $f_r^\mathbb N(s)=r$. 
By this and by Proposition \ref{prop:natural} there exists a homomorphism from $\mathbf A^\uparrow$ onto $(\mathbf B_{\bar r})^\uparrow$, so that $\theta_\mathbf A\subseteq \theta_{\mathbf B_{\bar r}}=\theta_{\hat r}$. In conclusion, $\theta_\mathbf A\subseteq \theta_\mathbf B=\bigcap_{r\in \mathsf b} \theta_{\hat r}$.


(2) $\Rightarrow$ (1) Since $\mathbf A^\uparrow\cong \mathcal T_\tau/\theta_\mathbf A$ and $\mathbf A\cong\mathbf A^\updownarrows$, then there exists an isomorphism  $g:(\mathcal T_\tau/\theta_\mathbf A)^\downarrow \cong \mathbf A$. Since $\mathcal T_\tau^\downarrow$ is generated by $\epsilon$, then  $\mathbf A$ is  generated by the thread $g^\mathbb N(\epsilon)=s$, so $\mathbf A =\mathbf A_{\bar s}$,  and $\theta_\mathbf A=\theta_{\hat s}$. 

Given $\mathbf B\in K$ of trace $\mathsf b$, we have that $\theta_\mathbf A \subseteq \theta_\mathbf B$, because $\mathbf A^\uparrow$ has  the initial property for $K^\smalltriangleup$.
For every $r\in \mathsf b$, we have to show that there exists a unique homomorphism $f$ from $\mathbf A$ into $\mathbf B_{\bar r}$ such that $f^\mathbb N(s)=r$.  We can define $f(t^\mathbf A(s))=t^\mathbf B(r)$, because
  $\theta_{\hat s}=\theta_\mathbf A \subseteq \theta_\mathbf B\subseteq  \theta_{\hat r}$. 
\end{proof}

The t-algebra $\mathcal T_{\tau}^\downarrow$ is the free t-algebra over $\epsilon$ in the class of all $\tau$-algebras.


\begin{theorem}\label{thm:ch2} Let $\mathcal D$ be a minimal  clone $\tau$-algebra and $K$ be the Et-variety generated by $\mathcal D^\downarrow$. Then the following conditions hold:
\begin{itemize}
\item[(i)] $\mathcal D$ has the initial property for $K^\smalltriangleup$.
\item[(ii)] $\mathcal D^\downarrow$ is the free $K$-t-algebra over $\epsilon^\mathcal D$.
\item[(iii)] If $\mathbf A\in K$, then $\mathbf A$ generates $K$, i.e., $K= \mathrm V_E(\mathbf A)$, iff $\mathcal D$ is isomorphic to $\mathbf A^\uparrow$.
\end{itemize} 
\end{theorem}

\begin{proof} By Proposition \ref{prop:thm} $\mathcal D^{\downuparrows} \cong \mathcal D$ and $\mathcal D^\downarrow\cong \mathcal D^{\downarrow\updownarrows}$.

(i) By $K=\mathrm V_E(\mathcal D^\downarrow)$ and Proposition \ref{prop:tre} we obtain that $\mathcal D^{\downuparrows}$ has the initial property for $K^\smalltriangleup$.  We get the conclusion from $\mathcal D^{\downuparrows} \cong \mathcal D$. 

(ii) By  $\mathcal D^\downarrow\cong \mathcal D^{\downarrow\updownarrows}$ and by Proposition \ref{prop:quattro} we get the conclusion.

(iii)
 ($\Rightarrow$)
By Proposition \ref{prop:tre} $\mathbf A^\uparrow\in K^\smalltriangleup$ has the initial property for $K^\smalltriangleup$. Since $\mathcal D^\downarrow\in K$, then $\mathcal D\cong \mathcal D^\downuparrows\in K^\smalltriangleup$. From  the minimality of $\mathcal D$  it follows that the unique homomorphism from $\mathbf A^\uparrow$ into $\mathcal D$ is an isomorphism, making $\mathcal D\cong\mathbf A^\uparrow$.
 ($\Leftarrow$)   By (i) $\mathbf A^\uparrow$ has the initial property for $K^\smalltriangleup$.
 From Proposition \ref{prop:tre} it follows that $K= \mathrm V_E(\mathbf A)$.
\end{proof}

\section{Topological Birkhoff and  Clone Algebras} \label{sec:ubandca}

The study of continuity properties of natural clone homomorphisms has been initiated  by Bodirsky and Pinsker \cite{BP15} for locally oligomorphic algebras.
Schneider \cite{Sc17} and Gehrke-Pinsker \cite{GP18} have recently shown that for any algebra $\mathbf B$ in the variety generated by an algebra $\mathbf A$, the induced natural clone homomorphism is uniformly continuous if and only if every finitely generated subalgebra of $\mathbf B$ is a homomorphic image of a subalgebra 
of a finite power of $\mathbf A$. In this section we generalise the topological Birkhoff's theorem to t-algebras. As a corollary, we get a new topological Birkhoff's theorem for algebras.

\subsection{Uniformity on $\mathsf{FCA}$s} Let $A$ be a set, $\mathsf a$ be a trace on $A$ and $A^{(\mathsf a)}=(A^\mathsf a,q_n^\mathsf a,\e_i^\mathsf a)$ be the pure full $\mathsf{FCA}$. 
We consider the uniformity $\mathsf U^A_\mathsf a$ on $A^\mathsf a$ generated by the equivalence relations
$$\alpha_s=\{ (\varphi,\psi):  \varphi(s)=\psi(s) \},\ \text{for every $s\in \mathsf a$}.$$

The uniformity $\mathsf U^A_\mathsf a$ is the least filter  in the powerset of $A^\mathsf a\times A^\mathsf a$ generated by the family of sets $\alpha_s$.
The uniformity $\mathsf U^A_\mathsf a$ satisfies the following condition, for every $\varphi,\psi, \gamma_1,\dots,\gamma_n\in A^\mathsf a$:
$$(\varphi,\psi)\in \alpha_{s[\gamma_1(s),\dots,\gamma_n(s)]}\ \text{iff}\ (q_n^\mathsf a(\varphi,\gamma_1,\dots,\gamma_n), q_n^\mathsf a(\psi,\gamma_1,\dots,\gamma_n)\in\alpha_s.$$


The uniformity $\mathsf U^A_\mathsf a$ defines the topology of pointwise convergence on $A^\mathsf a$ (see Section \ref{sec:us}). 
A subbase of the topology is given by the sets $\alpha_{\varphi,s}=\{ \psi:  \varphi(s)=\psi(s) \}$, for every $s\in \mathsf a$ and $\varphi\in A^\mathsf a$.



 Recall from Section \ref{sec:us}  the notation $\alpha[-]$ for an entourage $\alpha\in \mathsf U^A_\mathsf a$.

\subsection{A topological Birkhoff's theorem for t-algebras} \label{sec:tbtforta}
Let $A$ be a set and $\mathsf a$ be a trace on $A$. Then $\mathsf a^n$ is a trace on $A^n$, because, if $r= (r^1,\dots,r^n)\in\mathsf a^n$ then $r_i=(r^1_i,\dots,r^n_i)\in A^n$ for every $i\in\mathbb N$.
Moreover, if $a_i\in A^n$, then $r[a_1,\dots,a_k]=(r^1[a^1_1,\dots,a^1_k],\dots,r^n[a^n_1,\dots,a^n_k])\in \mathsf a^n$.

The following lemma is trivial.

\begin{lemma}
The full $\mathsf{FCA}$ $A^{(\mathsf a)}$ is isomorphic to a  $\mathsf{FCA}(A^n,\mathsf a^n)$ through  a mapping $\varphi:\mathsf a \to A \mapsto \varphi^n: \mathsf a^n\to A^n$, where $\varphi^n(r)=  (\varphi(r^1),\dots,\varphi(r^n))$
 for every $\varphi:\mathsf a\to A$ and $r=(r^1,\dots,r^n)\in \mathsf a^n$.
\end{lemma}

Let $\mathcal F$ be a pure $\mathsf{FCA}(A,\mathsf a)$.
We have two canonical  t-algebras of type $ F$:
\begin{enumerate}
\item $\mathbf A= (A,\mathsf a,\varphi)_{\varphi\in F}$.
\item $\mathbf F=(F,\pmb \epsilon^\mathbf F,\varphi^\mathbf F)_{\varphi\in F}$, where $\pmb \epsilon^\mathbf F=[\epsilon^\mathcal F]_\mathbb N$ is a basic trace and $$\varphi^\mathbf F(\epsilon^\mathcal F[\psi_1,\dots,\psi_k])=q_k^\mathsf a (\varphi,\psi_1,\dots,\psi_k)$$ for every $\epsilon^\mathcal F[\psi_1,\dots,\psi_k]\in \pmb \epsilon^\mathbf F$.
\end{enumerate}


\begin{lemma}\label{lem:r} For every $r=(r^1,\dots,r^n)\in \mathsf a^n$, the map $\bar r: \mathbf F\to \mathbf A^n$, defined by $\bar r(\varphi)=(\varphi(r^1),\dots,\varphi(r^n))$, is a t-homomorphism, whose kernel is the congruence $\alpha_{r^1}\cap\dots\cap \alpha_{r^n}\cap (F\times F)$ and whose image is the t-subalgebra  $(\mathbf A^n)_{\bar r}$ of  $\mathbf A^n$ generated by $r$. 
\end{lemma}

\begin{proof}
 By Proposition \ref{prop:zzz}(i) each map $\overline{r^i} : F\to A$ is a t-homomorphism from $\mathbf F$ into $\mathbf A$. By collecting all these t-homomorphisms we get the conclusion. 
\end{proof}


%




Let $\mathcal G$ be a pure $\mathsf{FCA}(B,\mathsf b)$ and $f: F\to G$ be an onto map such that $f(\e_i^\mathsf a)=\e_i^\mathsf b$. Then we have two other canonical  t-algebras of type $F$:
\begin{enumerate}
\item  $\mathbf B= (B,\mathsf b,f(\varphi))_{\varphi\in F}$.
\item $\mathbf G= (G, \pmb \epsilon^\mathbf G, f(\varphi)^\mathbf G)_{\varphi \in F}$, where $\pmb \epsilon^\mathbf G=[\epsilon^\mathcal G]_\mathbb N$ and, for all $u=\epsilon^\mathbf G[\psi_1,\dots,\psi_n]\in \pmb \epsilon^\mathbf G$,  $f(\varphi)^\mathbf G(u)=q_n^\mathsf b (f(\varphi),\psi_1,\dots,\psi_n)$.
\end{enumerate}



\begin{definition}
An onto homomorphism $f: \mathcal F\to \mathcal G$ of pure clone algebras is \emph{uniformly continuous} if, for every $s\in \mathsf b$  there exists  $r=(r^1,\dots,r^n)\in \mathsf a^n$ (for some $n\in\mathbb N$) such that
$$\forall \varphi,\psi \in F.\quad\bar r(\varphi)=\bar r(\psi) \Rightarrow \bar s(f\varphi)=\bar s(f\psi).$$
\end{definition}

\noindent In other words, $\varphi(r^1)=\psi(r^1) \land\dots\land \varphi(r^n)=\psi(r^n) \Rightarrow  (f\varphi)(s)= (f\psi)(s)$.

\begin{theorem}\label{thm:uchspfin} Let $\mathcal F$ be a  $\mathsf{FCA}( A,\mathsf a)$, $\mathcal G$ be a pure $\mathsf{FCA}(B,\mathsf b)$, $f:  F\to  G$ be an onto map such that $f(\e_i^\mathsf a)=\e_i^\mathsf b$, $\mathbf A= (A,\mathsf a,\varphi)_{\varphi\in F}$ and $\mathbf B= (B,\mathsf b, f(\varphi))_{\varphi\in F}$. Then the following conditions are equivalent:
\begin{enumerate}
\item $f$ is a uniformly continuous homomorphism from $\mathcal F$ onto $\mathcal G$.
\item For every $s\in \mathsf b$, there exist $r\in \mathsf a^n$ (for some $n$) and an onto t-homomorphism $h_s: (\mathbf A^n)_{\bar r}\to  \mathbf B_{\bar s}$  such that $h_s\circ \bar r= \bar s\circ f$.
\item Every t-subalgebra $\mathbf B_{\bar s}$ of $\mathbf B$ generated by $s\in\mathsf b$  is a t-homomorphic image of a  t-subalgebra of a finite t-power of $\mathbf A$.
\end{enumerate}
\end{theorem}

\begin{proof} (1) $\Rightarrow$ (2) We define $h_s(\bar r(\varphi))= \bar s(f\varphi)$ for every $\varphi\in F$. By uniform continuity $h_s$ is well defined.
Moreover, $h_s$ is a t-homomorphism, because $f:\mathbf F\to \mathbf G$, $\bar s: \mathbf G\to \mathbf B_{\bar s}$ and $\bar r:\mathbf F\to (\mathbf A^n)_{\bar r}$ (see Lemma \ref{lem:r}) are all onto t-homomorphisms.

(2) $\Rightarrow$ (3) Trivial.

(3)  $\Rightarrow$ (1) Let $h_s: \mathbf D\to \mathbf B_{\bar s}$ be an onto t-homomorphism of t-algebras of type $F$, where $\mathbf D$ is a t-subalgebra of $\mathbf A^n$ of trace $\mathsf d$. Since $h_s^\mathbb N: \mathsf d\to \mathsf b_s$ is onto, then there exists $r=(r^1,\dots,r^n)\in \mathsf d$ such that $h_s^\mathbb N(r)=s$. We claim that $(h_s)_{|(A^n)_{\bar r}}$, where $(A^n)_{\bar r}=\{\varphi(r): \varphi\in F\}$, is onto. Since $h_s$ is a t-homomorphism and $h_s\circ \varphi= (f\varphi)\circ h_s^\mathbb N$, then $$h_s(\varphi(r))= (f\varphi)(h_s^\mathbb N(r))=(f\varphi)(s).$$ It follows that $f$ is uniformly continuous.

We now prove that the map $f$ is a homomorphism of  clone algebras. Let $\psi_i(r)=(\psi_i(r^1),\dots,\psi_i(r^n))$. Then we have:
\[
\begin{array}{lll}
  &   & f(q_n^\mathsf a(\varphi,\psi_1,\dots,\psi_n))(s)  \\
    &  = &h_s(q_n^\mathsf a(\varphi,\psi_1,\dots,\psi_n)(r)) \\
  &  = &  h_s(\varphi(r^1[\psi_1(r^1),\dots,\psi_n(r^1)]),\dots,\varphi(r^n[\psi_1(r^n),\dots,\psi_n(r^n)])) \\
  &  = &   h_s(\varphi(r[\psi_1(r),\dots,\psi_n(r)]))\\
    &  = & (f\varphi)(h_s^\mathbb N(r[\psi_1(r),\dots,\psi_n(r)]))\\
     &  = & (f\varphi)(s[h_s(\psi_1(r)),\dots,h_s(\psi_n(r))])\\
      &  = & (f\varphi)(s[(f\psi_1)(s)),\dots,(f\psi_n)(s))])\\
      &  = &q_n^\mathsf b(f\varphi,f\psi_1,\dots,f\psi_n)(s).\qedhere
\end{array}
\]
\end{proof}

If $K$ is a class of t-algebras, then $\mathrm P_t^{\mathrm{fin}}(K)$ is the closure of $K$ by finite t-products.

\begin{theorem}\label{thm:tbta} {\rm (Topological Birkhoff for t-algebras)}
 Let $\mathbf A$ and $\mathbf  B$ be t-algebras of type $\tau$ and trace $\mathsf a$ and $\mathsf b$, respectively. The following
are equivalent.

(1) $\mathbf  B\in\mathrm V_E(\mathbf A)$ and the natural homomorphism from $\mathbf A^\uparrow$ onto $\mathbf B^\uparrow$  is uniformly continuous.

(2) Every t-subalgebra $\mathbf  B_{\bar s}$  of $\mathbf  B$ generated by $s\in \mathsf b$ is a t-homomorphic image of a  t-subalgebra of a finite t-power of $\mathbf A$.

\end{theorem}

\begin{proof} (1) $\Rightarrow$ (2) By Theorem \ref{thm:uchspfin} applied to $\mathcal F=\mathbf A^\uparrow$, $\mathcal G=\mathbf B^\uparrow$, $f:F\to G$ (defined by $t^\mathbf A\mapsto t^\mathbf B$), $(A,\mathsf a, t^\mathbf A)_{t^\mathbf A\in F}$ and $(B,\mathsf b, f(t^\mathbf A)=t^\mathbf B)_{t^\mathbf A\in F}$.
%
 
 (2) $\Rightarrow$ (1) 
%
%
 $\mathbf  B\in\mathrm V_E(\mathbf A)$ because $\mathrm V_E(\mathbf A)$ is closed under t-expansion. By Theorem \ref{thm:uchspfin} the natural homomorphism from $\mathbf A^\uparrow$ onto $\mathbf B^\uparrow$  is uniformly continuous.
\end{proof}


Note that $\mathbf  B_{\bar s}$ is generated by $s$ and $\mathrm{set}(s)$ can be countable infinite.

\subsection{A new topological Birkhoff's theorem for algebras}
We generalise the topological Birkhoff's theorem in \cite{Sc17} and \cite{GP18}.

\begin{definition}\label{def:t-extM}
  Let $\mathbf S$ be an algebra of finitary type $\rho$. The \emph{t-extension set of $\mathbf S$} is the set $\mathbf S^\star=\{\mathbf S_\mathsf a^\top:\ \text{$\mathsf a$ is a trace on $S$}\}$, where $\mathbf S_\mathsf a^\top$ is the  t-algebra  over $\mathbf S$ of trace $\mathsf a$ and type $\rho^\star$ defined in Example \ref{exa:extfinalg}.
\end{definition}

The t-algebras $\mathbf S_\mathsf a^\top$ of Definition \ref{def:t-extM} are $\rho$-dimensional and they are all full t-subalgebras of $\mathbf S_{S^\mathbb N}^\top$. The following lemma characterises the  t-extension set $\mathbf S^\star$ of $\mathbf S$.

\begin{lemma}\label{lem:i-v}  Let $\mathbf S$ be a $\rho$-algebra, $C=\mathrm{Clo}\mathbf S$ be the term clone  of $\mathbf S$ and $\mathbf A, \mathbf B\in \mathbf S^\star$. Then we have:


(i) The equational t-theory of $\mathbf A$ is $\theta_\mathbf A= \{t=u : t^\bullet = u^\bullet\in \mathrm{Th}_\mathbf S\}$.

(ii) 
$\mathbf A^\uparrow \cong \mathbf B^\uparrow\cong \mathcal C^\top_{\mathsf a,\rho^\star}$ for every trace $\mathsf a$ on $S$, 
where $\mathcal C^\top_{\mathsf a,\rho^\star}=(C^\top_\mathsf a,q_n^\mathsf a,\e_i^\mathsf a,f_\mathsf a^\top)_{f\in\rho^\star}$ is the block $\rho^\star$-algebra of trace $\mathsf a$ and universe $C^\top_\mathsf a=\{f_\mathsf a^\top : f\in C\}$ (see Definition \ref{def:blockalg}).

(iii) $\mathrm V_E(\mathbf A)=\mathrm V_E(\mathbf S^\star)$, for every $\mathbf A\in \mathbf S^\star$.
\end{lemma}

\begin{proof}
  
 (i) Since the t-algebras in $\mathbf S^\star$ are $\rho$-dimensional, then we can apply  Lemma \ref{lem:upterms2}(iv).
 
 (ii) By (i) $\theta_\mathbf A=\theta_\mathbf B$. Then $\mathbf A^\uparrow\cong \mathcal T_\tau/\theta_\mathbf A \cong \mathbf B^\uparrow$. The last part follows from Example \ref{lem:coincide}.
 
 (iii) By (i).
\end{proof}

%
%
%
%
%

All t-algebras in $\mathbf S^\star$ generate the same Et-variety. Moreover, if $K$ is an Et-variety, then either $\mathbf S^\star\subseteq K$ or $\mathbf S^\star\cap K=\emptyset$.

\begin{lemma}\label{lem:uc}
Let $\mathbf R$ and $\mathbf S$ be $\rho$-algebras, $\mathbf A\in\mathbf R^\star$ and $\mathbf B\in \mathbf S^\star$.
The natural homomorphism from $\mathbf A^\uparrow$ onto $\mathbf B^\uparrow$ exists and  is uniformly continuous iff the natural homomorphism from $\mathrm{Clo}\mathbf R$ onto $\mathrm{Clo}\mathbf S$ exists and is uniformly continuous.
\end{lemma}

\begin{proof} Let $\mathbf A=\mathbf R^\top_\mathsf a$ and $\mathbf B=\mathbf S^\top_\mathsf b$. Let $t,u\in T_{\rho^\star}(\e_1,\dots,\e_k)$, $r=(r^1,\dots,r^n)\in (R^\mathbb N)^n$ and $s\in S^\mathbb N$. Then by applying Lemma \ref{lem:upterms2} and the fact that $(-)^\bullet: T_{\rho^\star}\to F_\rho(I)$ is onto,  $\bar r(t^\mathbf A)=\bar r(u^\mathbf A) \Rightarrow  \bar s(t^\mathbf B)= \bar s(u^\mathbf B)$ iff $\bigwedge_{i=1}^n(t^\bullet)^\mathbf R(r^i_1,\dots,r^i_k)=(u^\bullet)^\mathbf R(r^i_1,\dots,r^i_k)$  $\Rightarrow  (t^\bullet)^\mathbf S(s_1,\dots,s_k)= (u^\bullet)^\mathbf S(s_1,\dots,s_k)$. For the converse, we can apply the same technique, because by Lemma \ref{lem:upterms1} every $\rho$-dimensional t-algebra satisfies   $t=(t^\bullet)^\star$.
\end{proof}

\begin{lemma}\label{lem:iff} Let 
$\mathbf R$  be a $\rho$-algebra  and $\mathbf R^\top_{R^\mathbb N}\in\mathbf R^\star$ be the t-algebra of complete trace $R^\mathbb N$. 
Then the following conditions hold:
\begin{enumerate}
\item  $(\mathbf R^n)^\top_{(R^n)^\mathbb N}= (\mathbf R^\top_{R^\mathbb N})^n$;
\item  $\mathbf S \leq \mathbf R$ iff $\mathbf S^\top_{S^\mathbb N}\leq_t  \mathbf R^\top_{R^\mathbb N}$;
\item  $f: \mathbf R \to \mathbf S$ is an onto homomorphism iff $f: \mathbf R^\top_{R^\mathbb N}\to \mathbf S^\top_{S^\mathbb N}$ is an onto t-homomorphism.
\end{enumerate}

\end{lemma}

\begin{proof} (1) By $(R^\mathbb N)^n= (R^n)^\mathbb N$.
The proof of (2) and (3) is trivial. 
\end{proof}

Let $\mathbf S$ be a $\rho$-algebra and $r\in S^\mathbb N$. We denote by $\langle r\rangle$ the subalgebra of $\mathbf S$ generated by $\mathrm{set}(r)$.  The subalgebra $\langle r\rangle$ is finitely generated if $|\mathrm{set}(r)|<\omega$.

\begin{theorem}\label{thm:gub} {\rm (Topological Birkhoff for algebras)} Let 
$\mathbf R$ and $\mathbf S$ be $\rho$-algebras, $\mathbf A\in\mathbf R^\star$ of complete trace $R^\mathbb N$, and $\mathbf B\in \mathbf S^\star$ of trace $\mathsf b$.
 The following are equivalent.
\begin{enumerate}
\item $\mathbf  B\in\mathrm V_E(\mathbf A)$ and the natural homomorphism from $\mathbf A^\uparrow$ onto $\mathbf B^\uparrow$  is uniformly continuous.
\item For every $r\in \mathsf b$, the t-subalgebra $\mathbf  B_{\bar r}$  of $\mathbf  B$ generated by $r$ belongs to $\mathrm H_t\mathrm S_t\mathrm P_t^{\mathrm{fin}}(\mathbf A)$.
\item $\mathbf S\in \mathrm{HSP}(\mathbf R)$ and the natural homomorphism from $\mathrm{Clo}\mathbf R$ onto $\mathrm{Clo}\mathbf S$  is uniformly continuous.
\item Every subalgebra $\langle r\rangle$ ($r\in \mathsf b$) of $\mathbf S$ belongs to $\mathrm{HSP}^{\mathrm{fin}}(\mathbf R)$.
\end{enumerate}
\end{theorem}

\begin{proof} The equivalence of (1) and (2) follows from Theorem \ref{thm:tbta}, while the equivalence of (1) and (3) follows from Lemma \ref{lem:uc}.
 
  (4) $\Leftrightarrow$ (2) 
By Lemma \ref{lem:upterms2} $B_{\bar r}=\{t^\mathbf B(r): t\in T_{\rho^\star}\} = \{(t^\bullet)^\mathbf S(s_1,\dots,s_n) : t\in T_{\rho^\star}, t^\bullet = t^\bullet(v_1,\dots,v_n)\ \text{for some}\ n\geq 0\}$. Then $\mathbf  B_{\bar r}$ and $\langle r\rangle$ have the same universe and $\mathbf  B_{\bar r}=\langle r\rangle^\top_{\mathsf b}$. Then the conclusion follows from Lemma \ref{lem:iff}.
  \end{proof}

The above theorem generalises all the known versions of Topological Birkhoff (see \cite{BP15,GP18, Sc17}), because we may have $|\mathrm{set}(r)|=\omega$ ($r\in \mathsf b$) in item (4) of Theorem \ref{thm:gub}.

We get the following corollary of Topological Birkhoff, where the equivalence of (1) and (2)  is due to Schneider \cite{Sc17} and Gehrke-Pinsker \cite{GP18}.

\begin{corollary} {\rm (Topological Birkhoff for algebras)} Let $\mathbf R$ and $\mathbf S$ be algebras of the same type and $s\in S^\mathbb N$. The following are equivalent.
\begin{enumerate}

\item $\mathbf S$ is contained in the variety generated by $\mathbf R$ and the natural homomorphism from $\mathrm{Clo}\mathbf R$ onto $\mathrm{Clo}\mathbf S$  is uniformly continuous.
\item Every finitely generated subalgebra of $\mathbf S$ is a homomorphic image of a subalgebra of a finite power of $\mathbf R$.
\item For every $a_1,\dots,a_n\in S$, the subalgebra  of $\mathbf S$ generated by $$\{a_1,\dots,a_n,s_{n+1},s_{n+2},\dots\}$$ is a homomorphic image of a subalgebra of a finite power of $\mathbf R$.
\end{enumerate}
If $S$ is an infinite set, then the above conditions are equivalent to
\begin{enumerate}
\item[(4)] Every countable infinitely generated subalgebra of $\mathbf S$ is a homomorphic image of a subalgebra of a finite power of $\mathbf R$.
\end{enumerate}
\end{corollary}

\begin{proof} It is sufficient to define  $\mathsf b$ as follows in Theorem \ref{thm:gub}: 
 (2)  $\mathsf b=\{r\in S^\mathbb N: |\mathrm{set}(r)| < \omega \}$; 
 (3)   $\mathsf b=[s]_\omega$.
 (4)   $\mathsf b=\{r\in S^\mathbb N: |\mathrm{set}(r)| = \omega \}$.
\end{proof}

%
%

\section{Conclusions}

All the original results presented in this paper stem from the definitions of clone algebra and t-algebra, which are, therefore, the main contributions of this work.  The results listed below give evidence of the relevance of the notion of clone algebra and t-algebra, which goes beyond providing a neat algebraic treatment of clones. 
\begin{itemize}
\item Theorems \ref{thm:ggg} and \ref{thm:giu} ensure that the varieties of clone algebras on the one hand and the (E,F)t-varieties on the other are strictly  connected.
 
\item In Theorems \ref{thm:bir1} and \ref{thm:bir2},  a Birkhoff-like Theorem has characterised those classes of t-algebras that are Et-varieties (resp. Ft-variety):
a class $K$ of t-algebras is an Et-variety (resp. Ft-variety) iff $K^\smalltriangleup$ is a variety of clone algebras and $K=K^{\smalltriangleup\blacktriangledown}$ (resp. $K=K^{\smalltriangleup\smalltriangledown}$).

\item Theorem \ref{thm:uchspfin} and other results in Section \ref{sec:tbtforta} have characterised  the uniformly continuous homomorphisms of functional clone algebras.

\item In Theorems \ref{thm:birkhoffforalgebras}  and \ref{thm:gub} we have given concrete examples that general results in universal clone algebra, when translated in terms of algebras and clones, give new versions of known theorems in universal algebra. We have applied this methodology to Birkhoff's HSP theorem and to the recent topological versions of Birkhoff's theorem.
\end{itemize}

The focus of the present paper is on the  relationship between clone algebras and t-algebras, and their meaning for  universal algebra and the theory of clones. A closer examination of potential implications  to universal algebra, the constraint satisfaction problem (see \cite{B15,BP16,BOP18}) and the  polymorphism clone of a structure (see \cite{BP15})
 is deferred to future work, which is currently in progress.

\end{document}